\documentclass[12pt]{article}

\usepackage{setspace, caption, subcaption, float, relsize, rotating, color, verbatim, multirow}
\usepackage{amsmath}
\usepackage{graphicx,psfrag,epsf}
\usepackage{enumerate}
\usepackage{natbib}
\usepackage{url}
\usepackage[colorlinks=true, urlcolor=blue,linkcolor=blue, citecolor=magenta, breaklinks=true]{hyperref}
\usepackage{empheq}
\usepackage{dsfont}
\usepackage{mathtools}
\usepackage{bm}
\usepackage{stmaryrd}
\usepackage{amssymb} 
\usepackage{amsfonts}
\usepackage{amsthm} 
\usepackage{mathrsfs}
\usepackage{bigints} 
\usepackage{enumerate} 
\usepackage{authblk} 
\usepackage{color} 
\usepackage{sectsty}
\usepackage{makecell}
\usepackage{amsfonts}
\usepackage{times}
\usepackage[plain,noend]{algorithm2e}

\newcommand{\E}{{\mathbb{E}}}

\newtheorem{thm}{Theorem}
\newtheorem{cor}{Corollary}
\newtheorem{prop}{Proposition}

\newtheorem{rk}{Remark}
\newtheorem{theorem}[thm]{Theorem}
\newtheorem{assumption}[]{Assumption}
\newtheorem{lemma}[thm]{Lemma}

\numberwithin{equation}{section}
\newcommand{\nc}{\newcommand}
\nc{\dps}{\displaystyle}
\nc{\tr}{\text{tr}}
\renewcommand{\thefootnote}{\fnsymbol{footnote}}

\def\squarebox#1{\hbox to #1{\hfill\vbox to #1{\vfill}}}
\def\boxit#1{\vbox{\hrule\hbox{\vrule\kern6pt
\vbox{\kern6pt#1\kern6pt}\kern6pt\vrule}\hrule}}

\newcommand{\blind}{0}

\def\boxit#1{\vbox{\hrule\hbox{\vrule\kern6pt
			\vbox{\kern6pt#1\kern6pt}\kern6pt\vrule}\hrule}}

\allowdisplaybreaks

\date{}

\begin{document}
\sloppy
\def\spacingset#1{\renewcommand{\baselinestretch}%
{#1}\small\normalsize} \spacingset{1}


\if0\blind
{
\title{{\bf Nonparametric Estimation of Mediation Effects with A General Treatment} \footnote{The authors are alphabetically ordered.}}
\author{Lukang Huang\thanks{E-mail: \texttt{lkhuang@nankai.edu.cn}}\\
School of Statistics \& Data Science, Nankai University\\
Wei Huang\thanks{E-mail: \texttt{wei.huang@unimelb.edu.au}}\\
School of Mathematics and Statistics, University of Melbourne\\
Oliver Linton\thanks{E-mail: \texttt{obl20@cam.ac.uk}}\\
Faculty of Economics, University of Cambridge\\
and \\
Zheng Zhang\thanks{\textit{Corresponding author}. E-mail: \texttt{zhengzhang@ruc.edu.cn}}\\
Center for Applied Statistics, Institute of Statistics \& Big Data, \\Renmin University of China
}
\maketitle
} \fi

\if1\blind
{
\bigskip
\bigskip
\bigskip
\medskip
} \fi

\bigskip
\begin{abstract}
To investigate causal mechanisms, causal mediation analysis decomposes the total treatment effect into the natural direct and indirect effects. This paper examines the estimation of the direct and indirect effects in a general treatment effect model, where the treatment can be binary, multi-valued, continuous, or a mixture. We propose generalized weighting estimators with weights estimated by solving an expanding set of equations. Under some sufficient conditions, we show that the proposed estimators are consistent and asymptotically normal. Specifically, when 
the treatment is discrete, the proposed estimators attain the semiparametric efficiency bounds. Meanwhile, when the treatment is continuous, the convergence rates of the proposed estimators are 
slower than $N^{-1/2}$; however, they are still more efficient than that constructed from the true weighting function. A simulation study reveals that our estimators exhibit a satisfactory finite-sample performance, while an application shows their practical value. 
\end{abstract}

\noindent%
\vfill

\renewcommand\thefootnote{\arabic{footnote}}
\newpage
\spacingset{1.45} 

\section{Introduction}
One essential goal of program evaluation and scientific research is to understand why and how a treatment variable affects the potential outcomes of interest, going beyond the estimation of the average treatment effects (ATE). 
For example, when assessing the effect of participating in an academic and vocational institute on criminal behavior, investigators might want to separate the program's direct effects on criminal behavior, such as integrity or discipline, from its effects relayed through the employment chances, which may indirectly affect the criminal behavior. In this regard, causal mediation analysis plays an important role by decomposing the total treatment effect into the natural direct effect and the indirect effects mediated through an intermediate variable, called the mediator \citep{robins1992identifiability,pearl2001direct}. Such an approach has been widely used in a number of disciplines in the medical and social sciences \citep[see e.g.,][]{baron1986moderator,imai:etal:11,vand:15}.

A fundamental problem of treatment effect analysis is that, although a randomized trial is the golden standard to identify the treatment effects, it is often unavailable or even unethical. More usually, in practice, treatments are randomly assigned to the individuals based on some features that also relate to the mediator and outcome of interest, thus causing the confounding issue. Such features are called confounders. The literature on causal mediation analysis, which can identify the direct and indirect effects from data with confounding, has rapidly grown over the last decade, producing abundant approaches and extensions.

An important nonparametric identifiability condition called the sequential ignorability assumption, which establishes a minimum set of assumptions required to identify the direct and indirect treatment effects regardless of the statistical models used, is commonly imposed \citep{imai2010general, imai:keel:yama:10,hsu2018nonparametric}. Most of the approaches focus on the binary treatment where an individual either receives or does not receive the treatment \citep[see, e.g.,][]{imai:keel:yama:10,joffe2007defining,vanderweele2009marginal,vanderweele2010odds,huber2014identifying,chan2016efficient,hsu2018nonparametric,liu2021large}. For example, based on the sequential ignorability assumption, \cite{hsu2018nonparametric} proposed nonparametric estimation of the binary natural direct and indirect effects by weighting observations using the inverse of estimated propensity scores. 

While binary treatment effects have been extensively studied, the estimation of continuous treatment effects has also drawn considerable attention (see \cite{hirano2004propensity, imai2004causal, galvao2015uniformly,kennedy2017non, Fong_Hazlett_Imai_2018,dong2019regression,colangelo2020double,Ai_Linton_Motegi_Zhang_cts_treat,huang2021unified} among others). However, all of the above studies focus on the total treatment effect rather than analyzing the causal mediation effects. In recent work, \cite{hsu2018direct} studied the identification and estimation of the natural direct and indirect effects when the treatment variable is continuous. Their estimators are constructed from weighting observations using the estimated marginal density of the treatment and the inverse of two estimated conditional densities, namely the conditional density of the treatment given the confounder and that given both the confounder and the mediator. \cite{singh2021kernel} identified the natural direct and indirect effects by the g-formula and estimated them using a Reproducing Kernel Hilbert Space approach. They showed their estimators are consistent and provided the rate of convergence, which achieves the minimax rate.

In many applications, the treatment variable can be complex, being a mixture of discrete and continuous elements, or even multidimensional, comprising both discrete and continuous components. For example, when evaluating the impact of an academic and vocational training program on criminal behavior, the causal effect might depend not only on participation in the program but also on the duration of this participation. Additionally, it could depend on multi-factors such as the duration and number of training courses undertaken.

The contribution of the work is three-fold. First, we propose a general framework for estimating the direct and indirect effects, which unifies the binary, multi-valued, and continuous treatments as well as the mixture of discrete and continuous treatments.   It can also be easily adapted to multidimensional treatment. Specifically, we integrate the sequential ignorability assumption for discrete treatments from \cite{hsu2018nonparametric} with that for continuous treatments in \cite{hsu2018direct}, and use a sieve method to estimate the treatment effect parameters. This synthesis results in a more general form of identifying and estimating direct and indirect effects, enabling nonparametric estimation for discrete, continuous, mixed and multidimensional treatments. 

Second, we demonstrate that the moment balancing method for estimating the stabilized weights, originally proposed for unconditional continuous treatment effects by \cite{Ai_Linton_Motegi_Zhang_cts_treat}, can improve the mediation analysis in the existing framework. Specifically, \citeauthor{hsu2018direct}'s (\citeyear{hsu2018direct}) mediation identification requires estimation of a weight that consists of the marginal density of the treatment and the inverse of two conditional densities, namely the conditional density of the treatment given the confounder and that given both the confounder and the mediator. \cite{hsu2018direct} estimate this weight using a ratio of kernel density estimators. However, it is known that the inverse probability weighting is sensitive to estimated densities in the denominator, potentially leading to extreme weights and unstable results (see e.g., \citealp{Fong_Hazlett_Imai_2018, Ai_Linton_Motegi_Zhang_cts_treat}). 
By reformulating the weight estimation, we apply the moment balancing idea to estimate the weight integrally, significantly reducing the likelihood of extreme weights and enhancing result accuracy.

Third, we verify that the proposed estimators of the direct and indirect effects are consistent and asymptotically converge to normal distributions after appropriate normalization. Specifically, the proposed estimators attain the semiparametric efficiency bounds when the treatment is discrete. Meanwhile, when the treatment is continuous, the convergence rates of the proposed nonparametric estimators are slower than $N^{-1/2}$. However, we show that our proposed estimators achieve an asymptotic variance smaller than or equal to those constructed from the true weighting function, provided that the tuning parameters are selected to ensure equivalent convergence rates. The equality is only attained when the conditional mean of the observed outcome given the observed treatment, mediators and confounders is always zero. 

The remainder of the paper is structured as follows. Section \ref{sec:basic_framework} establishes the basic
framework and Section \ref{sec:point estimation} presents the estimation
procedure. In Section \ref{sec:asymptotics}, we derive the asymptotic properties of the proposed estimators and their efficiency results. Section \ref{sec:variance} constructs the consistent estimators for the asymptotic variances based on plug-in approaches. In Section \ref{sec:tunning}, we propose data-driven approaches for selecting the tuning parameters. Subsequently, Section \ref{sec:simulation} reports the results of a simulation study and Section \ref{sec:application} applies the proposed
estimation method to analyze the effect of Job Corps, an educational and vocational training program, on
criminal activity. Finally, Section \ref{sec:conclusion} presents the study's conclusions.
\section {Basic Framework} \label{sec:basic_framework}
Let $T$ denote the observed treatment variable with support
$\mathcal{T}\subset\mathbb{R}$\footnote{Our causal mediation framework and the proposed estimation methods can also be adapted to the multiple dimensions of treatment. However, this article focuses on the univariate treatment variable for simplicity of presentation.}, where $\mathcal{T}$ is either a discrete
set, a continuum, or a mixture of discrete and continuum subsets. Let $f_T(t)$ denote a probability distribution function of $T$ at point $t$; that is, it stands for the probability density if $T$ is continuous at $t$ and the probability mass if $T$ is discrete at $t$. Under the standard framework of causal inference, let $\boldsymbol{M}(t)\in \mathcal{M} \subset \mathbb{R}^s$, for some positive integer $s$, denote a potential mediating variable that represents the value of the mediator if the treatment variable is equal to $t\in\mathcal{T}$. Similarly, let $Y(t,\boldsymbol{m}) \in \mathbb{R}$ denote the potential outcome if one receives treatment $t$ and mediator $\boldsymbol{m}$. The observed mediator and response are denoted by $\boldsymbol{M}:=\boldsymbol{M}(T)$ and $Y:=Y\{T,\boldsymbol{M}(T)\}$, respectively. 
Define $\mu(t,t'):=\mathbb{E}[Y\{t,\boldsymbol{M}(t')\}]$ for any $(t,t')\in\mathcal{T}\times\mathcal{T}$, then $\mu(t,t)-\mu(t',t')$ is the average treatment effect (ATE) caused by the change in the treatment level from $t'$ to $t$. In particular, if $\mathcal{T}=\{0,1\}$, then $\mu(1,1)-\mu(0,0)$ is the binary ATE studied by \cite{hahn1998role}, \cite{Hirano03}, \cite{chan2016globally}, \cite{hsu2018treatment}, and \cite{ai2022simple} among others; if $\mathcal{T}=\{0,1,...,J\}$ for some $J\in\mathbb{N}$, then $\mu(t,t)-\mu(t',t')$ is the multi-valued ATE studied by \cite{cattaneo2010efficient} and \cite{lee2018efficient} among others; and if the treatment variable is continuous, then $\mu(t,t)-\mu(t',t')$ is the continuous ATE studied by \cite{imai2004causal}, \cite{hirano2004propensity}, \cite{lee2010regression}, \cite{colangelo2020double}, and \cite{Ai_Linton_Motegi_Zhang_cts_treat} among others.

A primary goal of causal mediation analysis is to decompose the ATE into the average \emph{natural indirect effect} and the average \emph{natural direct effect}, for any $t\neq t'\in\mathcal{T}$,
\begin{align}\label{eq:ATEDecomposition}
\mu(t,t)-\mu(t',t')=\{\mu(t,t)-\mu(t',t)\}+\{\mu(t',t)-\mu(t',t')\},
\end{align}
where $\mu(t,t)-\mu(t',t)$ is the average natural direct effect representing the average difference if the treatment variable changes from $t'$ to $t$ while the mediator is held constant at $\boldsymbol{M}(t)$; $\mu(t',t)-\mu(t',t')$ is the average natural indirect effect representing the average difference if the mediator value changes from $\boldsymbol{M}(t')$ to $\boldsymbol{M}(t)$ while holding the treatment variable constant at $t'$. Therefore, this decomposition enables researchers to quantitatively explore the extent to which the treatment and mediator contribute to the treatment effect. 

Note from \eqref{eq:ATEDecomposition} that all the ATE, average natural indirect effect, and the average natural direct effect depend only on $\mu:\mathcal{T}\times \mathcal{T}\mapsto \mathbb{R}$. Our goal is then reduced to estimating $\mu$ on $\mathcal{T}\times \mathcal{T}$. 
Due to the confounding issue, the potential outcome $Y\{t,\boldsymbol{M}(t')\}$ and the potential mediator $\boldsymbol{M}(t)$ are not observed for all $(t,t')\in\mathcal{T}\times\mathcal{T}$. To address this identification problem, most studies impose a selection on the observable condition. Specifically, let $\boldsymbol{X}\in\mathcal{X}\subset\mathbb{R}^r$ be a vector of observed covariates, for some positive integer $r$. We maintain the following sequential ignorability assumption imposed on the treatment and mediator assignment \citep{imai:keel:yama:10,huber2014identifying,hsu2018nonparametric,hsu2018direct}.
\begin{assumption}
[\emph{Sequential Ignorability}]\label{as:SequentialIgnore} \
\begin{enumerate}[(i)]
\item $\{Y(t',\boldsymbol{m}), \boldsymbol{M}(t)\} \perp T|\boldsymbol{X}=\boldsymbol{x}$ for all $(t,t',\boldsymbol{m},\boldsymbol{x})\in \mathcal{T}\times \mathcal{T}\times \mathcal{M}\times \mathcal{X}$;
\item $Y(t',\boldsymbol{m}) \perp \boldsymbol{M}(t)|(T=t, \boldsymbol{X}=\boldsymbol{x})$ for all $(t,t',\boldsymbol{m},\boldsymbol{x})\in \mathcal{T}\times\mathcal{T}\times \mathcal{M}\times \mathcal{X}$,
\end{enumerate} 
where the conditional probability functions (in the sense as $f_T$) satisfy, $f_{T|\boldsymbol{X}}(t|\boldsymbol{x})>0$ and $f_{T|\boldsymbol{M},\boldsymbol{X}}(t|\boldsymbol{m},\boldsymbol{x})>0$ for all $(t,\boldsymbol{m},\boldsymbol{x})\in \mathcal{T}\times \mathcal{M}\times \mathcal{X}$.
\end{assumption}
When $\boldsymbol{M}$ is multidimensional, Assumption \ref{as:SequentialIgnore} needs to hold for each element in $\boldsymbol{M}$. Under Assumptions \ref{as:SequentialIgnore}, we show in the Appendix \ref{app:id} that $\mu(t,t')$ can be identified as follows:

\begin{align}\label{id:mutt'}
\mu(t,t')=&\mathbb{E}\left[\frac{\pi_{\boldsymbol{M},\boldsymbol{X}}(T,\boldsymbol{M},\boldsymbol{X})}{\pi_{\boldsymbol{M},\boldsymbol{X}}(T+\delta,\boldsymbol{M},\boldsymbol{X})} \cdot \pi_{\boldsymbol{X}}(T+\delta,\boldsymbol{X}) Y\bigg|T=t\right],
\end{align}
where, for $(t,t')\in \mathcal{T}\times\mathcal{T}$,
\begin{align*}
\delta:=t'-t, \ \pi_{\boldsymbol{Z}}(t,\boldsymbol{Z}):=\frac{f_T(t)}{f_{T|\boldsymbol{Z}}(t|\boldsymbol{Z})} \ \text{for} \ \boldsymbol{Z} \in \{\boldsymbol{X}, (\boldsymbol{M},\boldsymbol{X})\}\,.
\end{align*}
In the particular case of $\delta=0$, $\mu(t,t)=\mathbb{E}\left[\pi_{\boldsymbol{X}}(T,\boldsymbol{X}) Y|T=t\right]$ is the dose-response function studied in \cite{Ai_Linton_Motegi_Zhang_cts_treat}.

Let $\{T_i,\boldsymbol{M}_i,\boldsymbol{X}_i,Y_i\}_{i=1}^N$ denote an independent and identically distributed ($i.i.d.$) sample of observations drawn from the joint distribution of $(T,\boldsymbol{M},\boldsymbol{X},Y)$. If $\pi_{\boldsymbol{X}}(t,\boldsymbol{X})$ and $\pi_{\boldsymbol{M},\boldsymbol{X}}(t,\boldsymbol{M},\boldsymbol{X})$ were known, $\mu(t,t')$ can be estimated by the nonparametric series regression \citep{Newey94,Newey97}:
\begin{align}\label{def:muOracle}
\widetilde{\mu}(t,t^{\prime}):=&\left[\sum_{i=1}^{N} \frac{\pi_{\boldsymbol{M,X}}(T_i,\boldsymbol{M}_i,\boldsymbol{X}_i)}{\pi_{\boldsymbol{M,X}}(T_i+\delta,\boldsymbol{M}_i,\boldsymbol{X}_i)}\cdot \pi_{\boldsymbol{X}}(T_i+\delta,\boldsymbol{X}_i)Y_iu_{K_{0}}(T_i)^\top \right]\notag\\
&\times\left[ \sum_{i=1}^{N} u_{K_{0}}(T_i)u_{K_{0}}(T_i)^\top \right]^{-1}u_{K_{0}}(t),
\end{align}
where $u_{K_{0}}(T)=(u_{K_{0},1}(T),\ldots, u_{K_{0} ,K_{0}}(T))^{\top}$ is a prespecified basis function with dimension $K_{0}\in \mathbb{N}$. However, both $\pi_{\boldsymbol{X}}(T,\boldsymbol{X})$ and $\pi_{\boldsymbol{M},\boldsymbol{X}}(T,\boldsymbol{M},\boldsymbol{X})$ are unknown in practice and need to be replaced by some estimates.

\begin{rk} When $T$ is a continuous treatment variable, \cite{hsu2018direct} identify $\mu(t,t')$ as
		\begin{align}\label{Huber_mu_identify}
			\mu(t,t')=\mathbb{E}[Y(t,\boldsymbol{M}(t'))]=\lim_{h\to 0}\mathbb{E}\left[\frac{\mathcal{K}_h(T-t)}{f_{T|\boldsymbol{M},\boldsymbol{X}}(t|\boldsymbol{M},\boldsymbol{X})}\cdot \frac{f_{T|\boldsymbol{M,X}}(t'|\boldsymbol{M,X})}{f_{T|\boldsymbol{X}}(t'|\boldsymbol{X})}\cdot Y\right],
		\end{align}
		where $\mathcal{K}_h(x) = \mathcal{K}(x/h)/h$ with $\mathcal{K}(\cdot)$ a univariate kernel function. Our proposed identification differs in two significant ways. First, \eqref{id:mutt'} does not depend on the kernel weighting, $\mathcal{K}_h(T-t)$, allowing integration with any nonparametric regression methods for discrete, continuous or mixture $T$. Indeed, our proposed sieve estimator in \eqref{def:muOracle} accommodates all the scenarios. Second, we address the selection bias using three stabilized weights, $\pi_{\boldsymbol{M},\boldsymbol{X}}(T,\boldsymbol{M},\boldsymbol{X})$, $\pi_{\boldsymbol{M},\boldsymbol{X}}(T+\delta,\boldsymbol{M},\boldsymbol{X})$ and $\pi_{\boldsymbol{X}}(T+\delta,\boldsymbol{M},\boldsymbol{X})$, rather than the conditional densities used in \eqref{Huber_mu_identify}. This approach, leveraging the moment balancing idea from \cite{Ai_Linton_Motegi_Zhang_cts_treat} (see Section~\ref{sec:point estimation} for details), ensures stability, especially when the densities, $f_{T|\boldsymbol{M},\boldsymbol{X}}(t|\boldsymbol{M},\boldsymbol{X})$ and $f_{T|\boldsymbol{X}}(t'|\boldsymbol{X})$, are close to zero. However, \cite{Ai_Linton_Motegi_Zhang_cts_treat} does not consider the mediation and only requires one stabilized weight, $\pi_{\boldsymbol{X}}(T,\boldsymbol{X})$, to identify their causal parameter.
\end{rk}


\section{Estimation}\label{sec:point estimation}
In this section, we introduce a two-step nonparametric estimation of $\mu(t,t')$ across all pairs $(t,t')\in\mathcal{T}\times \mathcal{T}$ based on the identification \eqref{id:mutt'}. In the first step, shown in Section~\ref{sec:piEstimator}, we propose a unified framework for estimating the weighting functions, $\pi_{\boldsymbol{X}}$ and $\pi_{\boldsymbol{M},\boldsymbol{X}}$.  Subsequently, in the second step shown in Section~\ref{sec:muEstimator}, we regress the estimated $\big\{\pi_{\boldsymbol{M},\boldsymbol{X}}(T,\boldsymbol{M},\boldsymbol{X})\cdot \pi_{\boldsymbol{X}}(T+\delta,\boldsymbol{X})\cdot Y\big\} \big/ \pi_{\boldsymbol{M},\boldsymbol{X}}(T+\delta,\boldsymbol{M},\boldsymbol{X})$ against $T$ to obtain our final estimator of $\mu(t,t')$. For this estimation, we employ a sieve method suitable for discrete, continuous, or a mixture of discrete and continuous treatment types.

\subsection{Unified Framework for Estimating the Weighting Functions}\label{sec:piEstimator}
This section proposes the estimators for $\pi_{\boldsymbol{X}}$ and $\pi_{\boldsymbol{M},\boldsymbol{X}}$. A naive estimator of $\pi_{\boldsymbol{X}}$ can be the ratio of some estimated $f_T$ and $f_{T|\boldsymbol{X}}$. However, such a ratio estimator is very sensitive to small values of estimated $
f_{T|\boldsymbol{X}}$ and lead to undesirable results. To mitigate this problem, \cite{Ai_Linton_Motegi_Zhang_cts_treat} propose estimating $\pi_{\boldsymbol{X}}$ as a whole. 

We note that, using the definition of $\pi_{\boldsymbol{Z}}$ in \eqref{id:mutt'}, their approach can be applied to estimate $\pi_{\boldsymbol{X}}$ and $\pi_{\boldsymbol{M},\boldsymbol{X}}$ in the same framework. Note that for $\boldsymbol{Z}\in \{\boldsymbol{X}, (\boldsymbol{M},\boldsymbol{X})\}$,
\begin{align}
\mathbb{E}\left[ \pi_{\boldsymbol{Z}}(T,\boldsymbol{Z})u(T)v(\boldsymbol{Z})\right]
=\mathbb{E}[u(T)]\cdot\mathbb{E}[v(\boldsymbol{Z})], \label{moment1}
\end{align}
holds for any suitable functions $u(T)$ and $v(\boldsymbol{Z})$. Using \citet[Theorem 2]{Ai_Linton_Motegi_Zhang_cts_treat}, one can show that \eqref{moment1} identifies $\pi_{\boldsymbol{Z}}$. Equation~\eqref{moment1} suggests a possible way of estimating $\pi_{\boldsymbol{Z}}$; however, it implies an infinite number of equations, which is impossible to solve using a finite sample of observations. To overcome this difficulty, they approximate the infinite-dimensional function space using a sequence of finite-dimensional sieve spaces. Specifically, let $u_{k_{1}}(T)=(u_{k_{1},1}(T),\ldots$ $,u_{k_{1} ,k_{1}}(T))^{\top}$ and $v_{k_{\boldsymbol{Z}}}(\boldsymbol{Z})=\left( v_{k_{\boldsymbol{Z}}
,1}(\boldsymbol{Z}),\ldots,v_{k_{\boldsymbol{Z}},k_{\boldsymbol{Z}}}(\boldsymbol{Z})\right) ^{\top}$ be specified basis with dimensions $k_{1}\in\mathbb{N}$ and $k_{\boldsymbol{Z}}\in\mathbb{N}$, respectively, and let $K_{\boldsymbol{Z}}:=k_{1}\cdot k_{\boldsymbol{Z}}$. The 
functions $u_{k_{1}}(T)$ and $v_{k_{\boldsymbol{Z}}}(\boldsymbol{Z})$ are
approximation \emph{sieves} that can approximate any suitable functions $u(T)$ and $v(\boldsymbol{Z})$ arbitrarily well (see \cite{chen2007large} for discussions on the sieve approximation). 

Then \eqref{moment1} implies that
\begin{equation}
\mathbb{E}\left[ \pi_{\boldsymbol{Z}}(T,\boldsymbol{Z})u_{k_{1}}(T)v_{k_{\boldsymbol{Z}}}(\boldsymbol{Z})^{\top}\right] =\mathbb{E}[u_{k_{1}}(T)]\cdot \mathbb{E}[v_{k_{\boldsymbol{Z}}}(\boldsymbol{Z})]^{\top}. \label{sievemoment}
\end{equation}
Following \cite{Ai_Linton_Motegi_Zhang_cts_treat}, $\pi_{\boldsymbol{Z}}(T_i,\boldsymbol{Z}_i)$, for $i=1,\ldots,N$, can be estimated by the $\widehat{\pi}_i$'s, which are the following maximizer of an entropy, subject to the sample analog of \eqref{sievemoment}: 
{\small 
\begin{equation}
\left\{ 
\begin{array}{cc}
& \left\{ \widehat{\pi}_{i}\right\} _{i=1}^{N}=\arg\max\left\{ -N^{-1}
\sum_{i=1}^{N}\pi_{i}\log \pi_{i}\right\} \\[2mm] 
& \text{subject to}\ \frac{1}{N}\sum_{i=1}^{N}%
\pi_{i}u_{k_{1}}(T_{i})v_{k_{\boldsymbol{z}} }(\boldsymbol{Z}_{i})^{\top}=\left\{ \frac{1}{%
N}\sum_{i=1}^{N}u_{k_{1}} (T_{i})\right\} \left\{ \frac{1}{N}%
\sum_{j=1}^{N}v_{k_{\boldsymbol{z}}}(\boldsymbol{Z}_{j})^{\top}\right\} .%
\end{array}
\right. \label{E:cm1}
\end{equation}}
Note that by including a constant of one in
the sieve bases $u_{k_1}(T)$ and $v_{k_{\boldsymbol{Z}}}(\boldsymbol{Z})$, \eqref{E:cm1} implies that $N^{-1}\sum_{i=1}
^{N}\widehat{\pi}_{i}=1$. Moreover,
\begin{equation*}
\max\left( -N^{-1}\sum_{i=1}^{N}\pi_{i}\log \pi_{i}\right) =-\min\left\{ 
\sum_{i=1}^{N}(N^{-1}\pi_{i})\cdot\log\left(\frac{N^{-1}\pi_{i}}{N^{-1}}\right)\right\},
\end{equation*}
that is, the entropy maximization problem is equivalent to the minimization of the Kullback-Leibler divergence 
between $\{N^{-1}\pi_{i}\}_{i=1}^{N}$ and the uniform empirical distribution $\{N^{-1}\}_{i=1}^N$. Such an entropy serves a reasonable metric in the sense that the $\{N^{-1}\pi_i\}_{i=1}^N$ can be treated as a discrete probability distribution. This comes from the observation that $\pi_{\boldsymbol{Z}}(T_i,\boldsymbol{Z}_i)$ is positive and $\mathbb{E}\{\pi_{\boldsymbol{Z}}(T_i,\boldsymbol{Z}_i)\}=1$. Moreover, the formulation in \eqref{E:cm1} also guarantee the empirical counterparts, the $\pi_i$'s are all positive and $\sum^N_{i=1} (N^{-1} \pi_i) = 1$.

To estimate $\pi_{\boldsymbol{Z}}(T,\boldsymbol{Z})$, we use the dual solution to \eqref{E:cm1} showed by \cite{Ai_Linton_Motegi_Zhang_cts_treat}: 
\begin{align} \label{def:pihat_dual}
\widehat{\pi}_{K_{\boldsymbol{Z}}}(T,\boldsymbol{Z}):=\rho ^{\prime }\left\{
u_{k_{1}}(T)^{\top }\widehat{\Lambda}_{k_{1}\times k_{\boldsymbol{Z}}}v_{k_{\boldsymbol{Z}}}(\boldsymbol{Z})\right\} ,
\end{align}%
where $\rho'(v)=\exp(-v-1)$ is the first derivative of $\rho(v)=-\exp(-v-1)$, and $\widehat{\Lambda}_{k_{1}\times k_{\boldsymbol{Z}}}$ is the maximizer of the strictly
concave function $\widehat{G}_{k_{1}\times k_{\boldsymbol{Z}}}$ defined by
\begin{align}
&\widehat{G}_{k_{1}\times k_{\boldsymbol{Z}}}(\Lambda )\nonumber\\
&:=\frac{1}{N}\sum_{i=1}^{N}\rho \left\{u_{k_{1}}(T_{i})^{\top }\Lambda v_{k_{\boldsymbol{Z}}}(\boldsymbol{Z}_{i})\right\} -\left\{\frac{1}{N}\sum_{i=1}^{N}u_{k_{1}}(T_{i})\right\} ^{\top }\Lambda \left\{\frac{1}{N}\sum_{j=1}^{N}v_{k_{\boldsymbol{Z}}}(\boldsymbol{Z}_{j})\right\}.\label{def:G^hat}
\end{align}
The first order condition of \eqref{def:G^hat} implies that $\{\widehat{\pi}_{K_{\boldsymbol{Z}}}(T_{i},\boldsymbol{Z}_{i})\}_{i=1}^N$ satisfies the sample analog of \eqref{sievemoment}. Such restrictions improve the robustness of the estimation, with extreme weights being unlikely to be obtained. The concavity of \eqref{def:G^hat} enables us to easily obtain the solution via the Gauss-Newton algorithm. To ensure consistent estimation of $\pi_{\boldsymbol{Z}}(T,\boldsymbol{Z})$ for continuous $T$ (resp. $\boldsymbol{Z}$), the dimension of the bases, $k_1$ (reps. $k_{\boldsymbol{Z}}$), shall increases as the sample size increases. 


\subsection{Final Estimator of $\mu(t,t')$}\label{sec:muEstimator}
With \eqref{id:mutt'} and the estimated weighting functions $\widehat{\pi}_{K_{\boldsymbol{X}}}$ and $\widehat{\pi}_{K_{\boldsymbol{M,X}}}$, we define the estimator of $\mu(t,t')$ the sieve regression estimator:
\begin{align} \label{def:CBS}
\widehat{\mu}(t,t^{\prime}):=&\left[\sum_{i=1}^{N}\frac{ \widehat{\pi}_{K_{\boldsymbol{M,X}}}(T_i,\boldsymbol{M}_i,\boldsymbol{X}_i)}{\widehat{\pi}_{K_{\boldsymbol{M,X}}}(T_i+\delta,\boldsymbol{M}_i,\boldsymbol{X}_i)}\cdot \widehat{\pi}_{K_{\boldsymbol{X}}}(T_i+\delta,\boldsymbol{X}_i)Y_iu_{K_{0}}(T_i)^\top \right]\\
&\cdot\left[ \sum_{i=1}^{N} u_{K_{0}}(T_i)u_{K_{0}}(T_i)^\top \right]^{-1}u_{K_{0}}(t)\,. \notag
\end{align}

\begin{rk}
When $T$ is a binary treatment taking values in $\mathcal{T}=\{0,1\}$, $u_{k_1}(T)=u_{K_{0}}(T)=(\mathds{1}(T=0),\mathds{1}(T=1))^{\top}$ with $k_1=K_{0}\equiv 2$. In this case, $\{\mu(t,t'): (t,t')\in \mathcal{T}\times \mathcal{T}\}$ is a finite discrete set and estimable at a rate of $N^{-1/2}$. The semiparametric efficiency bounds and efficient estimation for $\{\mu(1,1),\mu(0,0)\}$ are presented in \cite{hahn1998role}, \cite{Hirano03}, and \cite{chan2016globally}. The semiparametric efficiency bounds and efficient estimation for the mediation causal effects $\{\mu(1,0),\mu(0,1)\}$ are presented in \cite{tchetgen2012semiparametric} and \cite{hsu2018nonparametric}. Specifically, \cite{hsu2018nonparametric} construct the estimator for $\mu(1,0)$ based on the following representation:
\begin{align*}
\mu(1,0)=\mathbb{E}[Y\{1,\boldsymbol{M}(0)\}]=\mathbb{E}\left[\frac{\mathds{1}(T=1)}{\mathbb{P}(T=1|\boldsymbol{M},\boldsymbol{X})}\frac{\mathbb{P}(T=0|\boldsymbol{M},\boldsymbol{X})}{\mathbb{P}(T=0|\boldsymbol{X})}\cdot Y\right],
\end{align*} 
with the generalized propensity score functions $\mathbb{P}(T=1|\boldsymbol{M},\boldsymbol{X})$ and $\mathbb{P}(T=1|\boldsymbol{X})$ estimated by the nonparametric series logit regression. \\
In the discrete setting, that is, $\mathcal{T}=\{0,1,...,J\}$ for some $J\in\mathbb{N}$, our proposed estimators can be obtained by taking $u_{k_1}(T)=u_{K_{0}}(T)=(\mathds{1}(T=0),\mathds{1}(T=1),\cdots,\mathds{1}(T=J))^{\top}$ with $k_1=K_{0}\equiv J+1$. We show in the next section that the proposed estimators have $\sqrt{N}$-asymptotic normality, provided some conditions on $k_{\boldsymbol{Z}}$, and attain the semiparametric efficiency bounds.\\
In a special case of a mixture of discrete and continuous treatment variable, where $T=0$ with positive probability and $T>0$ continuous with positive probability, our estimators can be obtained by taking $u_{k_1}(T)=u_{K_0}(T)=(\mathds{1}(T=0), u_{K_0,2}(T)\{1-\mathds{1}(T=0)\},\cdots, u_{K_0,K_0}(T)\{1-\mathds{1}(T=0)\})^{\top}$. \\
When $T$ is continuous or a mixture, we show in Theorem~\ref{thm:mutilde} that our proposed estimator converges to normal distributions slower than $N^{-1/2}$, but has an asymptotic variance smaller than or equal to that of $\widetilde{\mu}$ in \eqref{def:muOracle} that uses the true $\pi_{\boldsymbol{X}}$ and $\pi_{\boldsymbol{M,X}}$, provided the convergence rates are the same. The equivalence is only taken when $\mathbb{E}(Y|T,\boldsymbol{M},\boldsymbol{X})=0$.
\end{rk}

\begin{rk} For continuous $T$, \cite{hsu2018direct} propose a nonparametric weighted kernel (Nadaraya-Watson) type estimator for $\mu(t,t')$ based on \eqref{Huber_mu_identify}. The conditional density functions $f_{T|\boldsymbol{M},\boldsymbol{X}}$ and $f_{T|\boldsymbol{X}}$ are estimated through the kernel method. Their estimators are shown to be asymptotically normal.\\
In such a case, using our proposed estimators of the weights $\pi_{\boldsymbol{Z}}$, a kernel (Nadaraya-Watson) type estimator for $\mu(t,t')$, an alternative to \eqref{def:CBS}, can also be constructed: 
\begin{align}\label{def:CBK}
\widehat{\mu}_h(t,t')=\frac{\sum_{i=1}^N\widehat{R}_{K_{\boldsymbol{M,X}}}(T_i,T_i+\delta,\boldsymbol{M}_i,\boldsymbol{X}_i) \cdot \widehat{\pi}_{K_{\boldsymbol{X}}}(T_i+\delta,\boldsymbol{X}_i) \mathcal{K}_h(T_i-t)Y_i }{\sum_{i=1}^N\widehat{R}_{K_{\boldsymbol{M,X}}}(T_i,T_i+\delta,\boldsymbol{M}_i,\boldsymbol{X}_i) \cdot \widehat{\pi}_{K_{\boldsymbol{X}}}(T_i+\delta,\boldsymbol{X}_i) \mathcal{K}_h(T_i-t) }\,,
\end{align}
where $\widehat{R}_{K_{\boldsymbol{M,X}}}(t,t',\boldsymbol{m},\boldsymbol{x}):= \widehat{\pi}_{K_{\boldsymbol{M,X}}}(t,\boldsymbol{m},\boldsymbol{x})/\widehat{\pi}_{K_{\boldsymbol{M,X}}}(t',\boldsymbol{m},\boldsymbol{x})$ for $(t,t',\boldsymbol{m},\boldsymbol{x})\in \mathcal{T}\times\mathcal{T}\times\mathcal{M}\times\mathcal{X}$.
We establish its asymptotic properties in Appendix \ref{app:asymptoics_kernel}, where we show that it attains a convergence rate as a standard Nadaraya-Watson estimator and it  achieves a faster convergence rate and smaller asymptotic variance than \cite{hsu2018direct}, unless \cite{hsu2018direct} use a same bandwidth for all the estimations of $f_{T|\boldsymbol{M},\boldsymbol{X}}$, $f_{T|\boldsymbol{X}}$ and $\mu(t,t')$. However, in practice, since these three quantities are conditional on different set of variables, namely $\{\boldsymbol{M},\boldsymbol{X}\}$, $\boldsymbol{X}$ and $T$, respectively, using a same bandwidth does not guarantee good finite sample performance. Indeed, \cite{hsu2018direct} suggest using different bandwidths for these estimations, and thus their $\mu(t,t')$ estimator cannot achieve the optimal convergence rate (see footnote 7 in \citealp{hsu2018direct}). Our proposed estimator can achieve a faster rate and smaller variance, because our proposed $\widehat{\pi}_{\boldsymbol{Z}}$ efficiently incorporates the additional information available in the form of $\pi_{\boldsymbol{Z}}$ as per \eqref{moment1} (see more discussions in \citealp{Hirano03}).

\end{rk}

\begin{rk}In addition to the natural direct effect $\mathbb{E}[Y(t,\boldsymbol{M}(t))]-\mathbb{E}[Y(t',\boldsymbol{M}(t))]$ for $t\neq t'$ studied in this paper, the controlled direct effect (CDE) defined as $\mathbb{E}[Y(t,\boldsymbol{m})]-\mathbb{E}[Y(t',\boldsymbol{m})]$ for $t\neq t'$ and a fixed  value $\boldsymbol{m}$ of the mediator is also of interest in the literature, see
	\cite{goetgeluk2008estimation}, \cite{vanderweele2009marginal}, \cite{hong2015ratio} for examples. The proposed method in this paper can also be adapted to estimate $\mathbb{E}[Y(t,\boldsymbol{m})]$ as well as CDE; indeed, by \citet[Theorem 2]{hong2015ratio},  $\mathbb{E}[Y(t,\boldsymbol{m})]$ can be identified as the following density ratio weighting form:
	\begin{align*}
		\mathbb{E}[Y(t,\boldsymbol{m})]=\mathbb{E}\left[\frac{f_{T,\boldsymbol{M}}(T,\boldsymbol{M})}{f_{T,\boldsymbol{M}|\boldsymbol{X}}(T,\boldsymbol{M}|\boldsymbol{X})}Y\bigg|T=t,\boldsymbol{M}=\boldsymbol{m}\right], 
	\end{align*}
which is similar to the identification of $\mathbb{E}[Y(t,\boldsymbol{M}(t'))]$ in Eq. \eqref{id:mutt'}. Therefore, we can apply the same procedure of estimating $\pi_{\boldsymbol{Z}}(T,\boldsymbol{Z})$ and $\mathbb{E}[Y(t,\boldsymbol{M}(t'))]$ to obtain the estimators of $f_{T,\boldsymbol{M}}(T,\boldsymbol{M})/f_{T,\boldsymbol{M}|\boldsymbol{X}}(T,\boldsymbol{M}|\boldsymbol{X})$ and $\mathbb{E}[Y(t,\boldsymbol{m})]$, respectively.
\end{rk}

\section{Large Sample Properties} \label{sec:asymptotics}
This section studies the asymptotic properties of the proposed estimator $\widehat{\mu}(t,t')$. The convergence rates for $\widehat{\pi}_{K_{\boldsymbol{Z}}}(\cdot ,\boldsymbol{Z})$, $\boldsymbol{Z}\in\{\boldsymbol{X},(\boldsymbol{M},\boldsymbol{X})\}$, are implied directly by the results in \cite{Ai_Linton_Zhang_CDQ}. We also recall these results in Appendix \ref{appendix:preliminary}. 
To facilitate the presentation our main results, we
introduce the following notations: $\Phi_{K_{0}\times K_{0}}:=\mathbb{E}[u_{K_{0}}(T)u_{K_{0}}^{\top}(T)]$, 
\begin{align}
d_{K_{0},i}(T,\boldsymbol{M},\boldsymbol{X},Y;\delta): =&\text{IF}_{\pi_{\boldsymbol{X}},i}\bigg\{\delta,\frac{\pi_{\boldsymbol{M},\boldsymbol{X}}(T,\boldsymbol{M},\boldsymbol{X})}{\pi_{\boldsymbol{M},\boldsymbol{X}}(T+\delta,\boldsymbol{M},\boldsymbol{X})}\cdot u_{K_0}(T) Y \bigg\} \label{def:dki}\\
&+\text{IF}_{\pi_{\boldsymbol{M},\boldsymbol{X}},i}\bigg\{0,\frac{\pi_{\boldsymbol{X}}(T+\delta,\boldsymbol{X})}{\pi_{\boldsymbol{M},\boldsymbol{X}}(T+\delta,\boldsymbol{M},\boldsymbol{X})}\cdot u_{K_0}(T) Y \bigg\} \notag\\
&-\text{IF}_{\pi_{\boldsymbol{M},\boldsymbol{X}},i}\bigg\{\delta,\frac{\pi_{\boldsymbol{M},\boldsymbol{X}}(T,\boldsymbol{M},\boldsymbol{X})\pi_{\boldsymbol{X}}(T+\delta,\boldsymbol{X})}{\pi^2_{\boldsymbol{M,X}}(T+\delta, \boldsymbol{M,X})}\cdot u_{K_0}(T) Y\bigg\} \notag\\
&-\mathbb{E}\left[\frac{\pi_{\boldsymbol{M},\boldsymbol{X}}(T_i,\boldsymbol{M}_i,\boldsymbol{X}_i)\pi_{\boldsymbol{X}}(T_i+\delta,\boldsymbol{X}_i)}{\pi_{\boldsymbol{M,X}}(T_i+\delta, \boldsymbol{M_i,X_i})}\cdot u_{K_0}(T_i) Y_i\Big|T_i\right]\notag\\
&+\mathbb{E}\left[\frac{\pi_{\boldsymbol{M},\boldsymbol{X}}(T_i,\boldsymbol{M}_i,\boldsymbol{X}_i)\pi_{\boldsymbol{X}}(T_i+\delta,\boldsymbol{X}_i)}{\pi_{\boldsymbol{M,X}}(T_i+\delta, \boldsymbol{M_i,X_i})}\cdot u_{K_0}(T_i) Y_i\right],\notag
\end{align}
for any $\delta\geq 0$ and $i=1,\ldots, N$, where the $\text{IF}_{\pi_{\boldsymbol{Z}},i}$'s are the i.i.d. mean zero influence functions for $\widehat{\pi}_{\boldsymbol{Z}}$, $\boldsymbol{Z}\in \{\boldsymbol{X},(\boldsymbol{M},\boldsymbol{X})\}$, defined in \eqref{def:IF_piZ} in Appendix~\ref{appendix:preliminary}, and
\begin{align}
V_{tt'}:= &\mathbb{E}\left[ \left\{ u_{K_{0}}^{\top}(t)\Phi_{K_{0}\times K_{0}}^{-1}d_{K_{0},i}(T,\boldsymbol{M},\boldsymbol{X},Y;\delta)\right\}^{2}\right] \label{Vtt'}\\
= & u_{K_{0}}^{\top}(t)\cdot\Phi_{K_{0}\times K_{0}}^{-1}\cdot\mathbb{E}\left[ d_{K_{0},i}(T,\boldsymbol{M},\boldsymbol{X},Y;\delta)d_{K_{0}, i}^{\top}(T,\boldsymbol{M},\boldsymbol{X},Y;\delta)\right] \cdot\Phi_{K_{0}\times K_{0}}^{-1}\cdot u_{K_{0}}(t). \notag
\end{align}
The following conditions are maintained throughout this article. 
\begin{assumption}\label{as:appprox_mu} For every fixed $t'\in\mathcal{T}$, there exist a $\gamma^{*}\in\mathbb{R}^{K_{0}}$ and a positive constant $\beta>0$ such that
\[\sup_{t\in\mathcal{T} }|\mu(t,t')-(\gamma^{*})^{
\top}u_{K_{0}}(t)| =O\left(K_{0}^{-\beta}\right).
\]
\end{assumption} 
\begin{assumption}\label{as:ddmatrix} The eigenvalues of $\mathbb{E}\left[u_{K_{0}}(T)u_{K_{0}}^{\top
}(T)\right]$ and $\mathbb{E}\left[ d_{K_{0},i}(T,\boldsymbol{M},\boldsymbol{X}, Y ;\delta)\cdot \notag\right. \\ \left. d^{\top}_{K_{0},i}(T,\boldsymbol{M},\boldsymbol{X},Y;\delta )\right]$ are bounded away from zero and infinity uniformly with respect to $K_0\in\mathbb{N}$.
\end{assumption}
Assumption \ref{as:appprox_mu} requires the sieve approximation error of the function
$\mu(\cdot,t')$ to shrink at a
polynomial rate. This condition is satisfied for a variety of sieve basis
functions. For example, if $T$ is discrete, then
the approximation error is zero for sufficiently large $K_0$; thus,
in this case, Assumption \ref{as:appprox_mu} is satisfied with $\beta=+\infty$. If $T$ is continuous or a mixture, the polynomial rate $\beta$ depends positively on the smoothness of $\mu(t,t')$ in $t$; indeed, for power series and $B$-splines, $\beta$ is the smoothness of $\mu(t;t')$ in $t$ \citep[][Section~2.3.1]{chen2007large}.
We show that the convergence rate of the estimated $\mu(t,t')$ is bounded by this polynomial rate. Assumption \ref{as:ddmatrix} essentially ensures the variance of the
estimator is non-degenerate.
Under these conditions and Assumptions \ref{as:suppX} -- \ref{as:u&v} presented in Appendix \ref{appendix:preliminary}, we establish the following two theorems, which hold whenever the treatment is continuous, discrete or mixed:
\begin{theorem}
\label{thm:sieve} Under Assumptions
\ref{as:SequentialIgnore} -- \ref{as:ddmatrix} and \ref{as:suppX} -- \ref{as:u&v} presented in Appendix \ref{appendix:preliminary},we have
\begin{enumerate}
\item (Convergence Rates) 
\begin{align*}
&\int_{\mathcal{T}}|\widehat{\mu}(t,t')-\mu(t,t')|^2dF_{T}(t)
=O_p\Bigg(\left\{K_0^{-2\beta}+\frac{K_0}{N}\right\}+\zeta(K_{\boldsymbol{X}})^2\left\{\frac{K_{\boldsymbol{X}}}{N}+K_{\boldsymbol{X}}^{-2\alpha_{\boldsymbol{X}}}\right\}\\
&\qquad \qquad \qquad \qquad \qquad \qquad \qquad \qquad +\zeta(K_{\boldsymbol{M,X}})^2\left\{\frac{K_{\boldsymbol{M,X}}}{N}+K_{\boldsymbol{M,X}}^{-2\alpha_{\boldsymbol{M,X}}}\right\}\Bigg),\\
&\sup_{t\in\mathcal{T}}|\widehat{\mu}(t,t')-\mu(t,t')|=O_p\Bigg(\zeta_1(K_0)\Bigg\{\left(K_0^{-\beta}+\frac{K_0}{N}\right)+\zeta(K_{\boldsymbol{X}})\left(\sqrt{\frac{K_{\boldsymbol{X}}}{N}}+K_{\boldsymbol{X}}^{-\alpha_{\boldsymbol{X}}}\right)\\
&\qquad \qquad \qquad \qquad \qquad \qquad \qquad \qquad +\zeta(K_{\boldsymbol{M,X}})\left(\sqrt{\frac{K_{\boldsymbol{M,X}}}{N}}+K_{\boldsymbol{M,X}}^{-\alpha_{\boldsymbol{M,X}}}\right)\Bigg\}\Bigg),
\end{align*}
hold for every fixed $t'\in\mathcal{T}$, where $\zeta(K_{\boldsymbol{Z}})$, for $\boldsymbol{Z}=\{\boldsymbol{X},(\boldsymbol{M},\boldsymbol{X})\}$, and $\zeta_1(K_0)$ are defined in Assumption \ref{as:u&v} of Appendix~\ref{appendix:preliminary}. 
\item (Asymptotic Normality) Suppose $\sqrt
{N}K_{0}^{-\beta}\to 0$. Then, for any fixed $t,t'\in\mathcal{T}$,
\[
\sqrt{N}\left\{ \widehat{\mu}(t,t')-\mu(t,t') \right\}=\frac{1}{\sqrt{N}}\sum_{i=1}^{N}\phi_{tt'}(Y_i,T_i,\boldsymbol{M}_i,\boldsymbol{X}_i;\delta)+o_p(1),
\]
where 
\[\phi_{tt'}(Y_i,T_i,\boldsymbol{M}_i,\boldsymbol{X}_i;\delta)=u_{K_0}^{\top}(t)\Phi_{K_0\times K_0}^{-1}d_{K_0,i}(T,\boldsymbol{M},\boldsymbol{X},Y;\delta)\,.\]
Thus, we have $$\sqrt{N}V_{tt'}^{-1/2}\{ \widehat{\mu}(t,t')-\mu(t,t')\}\xrightarrow{d}N(0,1)\,,$$
and $V_{tt'}= \text{const}\times \|u_{K_0}(t)\|^2= O(K_0)$ for some positive constant $\text{const}$. 
\end{enumerate}
\end{theorem}
The proof of Theorem \ref{thm:sieve} is presented in section~S3 in the supplemental material. In addition to the convergence rate and asymptotic normality, we also provide the asymptotic linear
expansion of $\sqrt{N}\{ \widehat{\mu}(t,t')-\mu(t,t')\}$. This can help conduct statistical inference, as we can approximate the limiting distribution of our estimator by adopting the exchangeable
bootstrap method \citep{chernozhukov2013inference,Donald2014Estimation,huang2021unified}. Using Theorem \ref{thm:sieve}, the asymptotic results of the estimated direct and indirect effects can be easily obtained. In particular, note from Remark~1 and Assumption~\ref{as:appprox_mu} that if the treatment is discrete, $K_0$ is a constant independent of $N$ and $\beta = +\infty$; thus $\widehat{\mu}(t,t')-\mu(t,t')$ attains a $\sqrt{N}$-asymptotic normality. If the treatment is continuous or a mixture, and $u_{K_0}(t)$ is a sieve series, then the convergence rate of $\widehat{\mu}(t,t')-\mu(t,t')$ is slower than $N^{-1/2}$; see Remark~\ref{rk:CurseofDimensionality} for more discussion on the curse of dimensionality that arises from $\boldsymbol{M}$ and $\boldsymbol{X}$.
The next theorem shows that our proposed estimator of mean potential outcome $\widehat{\mu}(t,t')$ is asymptotically more efficient than the oracle estimator $\widetilde{\mu}(t,t')$ in \eqref{def:muOracle} that uses the true $\pi_{\boldsymbol{X}}$ and $\pi_{\boldsymbol{M,X}}$.

\begin{theorem}\label{thm:mutilde}
Suppose $\sqrt{N}K_0^{-\beta}\to 0$. Under Assumptions~\ref{as:SequentialIgnore} -- \ref{as:ddmatrix}, for any fixed $t,t'\in \mathcal{T}$,
$$\sqrt{N}\widetilde{V}_{tt'}^{-1/2}\{ \widetilde{\mu}(t,t')-\mu(t,t')\}\xrightarrow{d}N(0,1)\,,$$
where $\widetilde{V}_{tt'}$ is the asymptotic variance satisfying
$\widetilde{V}_{tt'}\ge V_{tt'}\,.$
\end{theorem}
It is known that in the estimation of average treatment effects with a binary treatment, the inverse probability weighting (IPW) estimator constructed from nonparametrically estimated propensity score is more efficient than that constructed by using the true one, see \cite{hahn1998role}, \cite{Hirano03} and \cite{chen2008semiparametric}. Theorem \ref{thm:mutilde} establishes the similar result for the mediation effects with a general treatment. The proof of Theorem~\ref{thm:mutilde} is presented in section~S4 in the supplemental material, where we also derive the asymptotic linear expansion of $\sqrt{N}\widetilde{V}_{tt'}^{-1/2}\{ \widetilde{\mu}(t,t')-\mu(t,t')\}$ and the detailed asymptotic variance $\widetilde{V}_{tt'}$. Theorem~\ref{thm:mutilde} implies the efficiency of our estimator over the oracle one for discrete, continuous, and mixed treatments. Specifically, when the treatment variable is discrete, we further prove the corollary below that our proposed estimators attain the semiparametric efficiency bounds. The proof is in section~S5 in the supplemental material.
\begin{cor}\label{cor:efficiency}
Suppose that the treatment variable $T$ is discrete with values in $\mathcal{T}=\{0,1,...,J\}$, where $J\geq 1$ is an positive integer. Under Assumptions
\ref{as:SequentialIgnore} -- \ref{as:ddmatrix} and \ref{as:suppX} -- \ref{as:u&v} presented in Appendix \ref{appendix:preliminary}, we have
$$\sqrt{N}(\widehat{\mu}(t,t')-\mu(t,t'))\stackrel {d} {\longrightarrow} \mathcal{N}(0,
V_{tt'}),$$ where $V_{tt'}=\mathbb{E}\left[S^2_{\mu(t,t')}\right]$ and
\begin{eqnarray}
S_{\mu(t,t')} & = &\frac{\mathds{1}\{T=t\}f_{\boldsymbol{M}\mid T,\boldsymbol{X}}(\boldsymbol{M}\mid
T=t',\boldsymbol{X})}{f_{T\mid \boldsymbol{X}}(t\mid \boldsymbol{X})f_{\boldsymbol{M}\mid T,\boldsymbol{X}}(\boldsymbol{M}\mid T=t,\boldsymbol{X})}\{Y-\E(Y\mid \boldsymbol{X},\boldsymbol{M},T=t)\}\nonumber \\
&&+\frac{\mathds{1}\{T=t'\}}{f_{T\mid \boldsymbol{X}}(t'\mid
\boldsymbol{X})}\{\E(Y\mid \boldsymbol{X},\boldsymbol{M},T=t)-\eta(t,t',\boldsymbol{X})\}\notag\\
&&+\eta(t,t',\boldsymbol{X})-\mu(t,t'), \notag
\end{eqnarray}
where
\begin{align*}
\eta(t,t',\boldsymbol{X}) = \int \E(Y\mid \boldsymbol{X},\boldsymbol{M}=\boldsymbol{m},T=t)f_{\boldsymbol{M}\mid T,\boldsymbol{X}}(\boldsymbol{m}\mid T=t',\boldsymbol{X})d\boldsymbol{m}\,,
\end{align*}
and $S_{\mu(t,t')}$ is equivalent to the efficient influence function in \cite{hsu2018nonparametric}.
\end{cor}

\begin{rk}\label{rk:CurseofDimensionality}
Note from Assumptions~\ref{as:smooth_pi} and \ref{as:u&v} in Appendix~\ref{appendix:preliminary} that the more continuous or mixed components in the covariate set $\boldsymbol{X}$ or the mediator set $\boldsymbol{M}$, the harder our estimation of $\pi_{\boldsymbol{Z}}$ would be, for $\boldsymbol{Z}\in\{\boldsymbol{X},(\boldsymbol{M},\boldsymbol{X})\}$. This is a curse of dimensionality of our fully nonparametric estimation. To resolve this issue, it is common to reduce the dimension. A solution is to use an index model \cite[see e.g.][]{fan1996local, lewbel2007nonparametric}. Specifically, we may assume that the random treatment assignment $T$ depends on $\boldsymbol{Z}$ through a function known up to a finite number of parameters $\boldsymbol{\theta}$, $s_{\boldsymbol{Z},\boldsymbol{\theta}}:\mathcal{Z}\mapsto \mathbb{R}$, in the sense that $f_{T|\boldsymbol{Z}}(t|\boldsymbol{z}) = f_{T|s_{\boldsymbol{Z},\boldsymbol{\theta}}}\{t|s_{\boldsymbol{Z},\boldsymbol{\theta}}(\boldsymbol{z})\}$ for all $(t,\boldsymbol{z})\in \mathcal{T}\times\mathcal{Z}$, and for $\boldsymbol{Z}\in\{\boldsymbol{X},(\boldsymbol{M},\boldsymbol{X})\}$ that has a large dimension of continuous or mixed components. Then $\pi_{\boldsymbol{Z}}$ can be identified by 
\begin{equation}\label{eq:idIndex}
\mathbb{E}[\pi_{\boldsymbol{Z}}(T,\boldsymbol{Z})u(T)v\{s_{\boldsymbol{Z},\boldsymbol{\theta}}(\boldsymbol{Z})\}] = \mathbb{E}[u(T)]\cdot \mathbb{E}[v\{s_{\boldsymbol{Z},\boldsymbol{\theta}}(\boldsymbol{Z})\}]
\end{equation}
for all suitable functions $u$ and $v$. However, the moment equation does not identify $\boldsymbol{\theta}$, we need to first estimate $\boldsymbol{\theta}$. One possible way is to estimate $\boldsymbol{\theta}$ by the maximum likelihood (MLE) with a nonparametric estimator of $f_{T|s_{\boldsymbol{Z},\boldsymbol{\theta}}}(t|s)$, $\widehat{f}_{T|s_{\boldsymbol{Z},\boldsymbol{\theta}}}(t|s)$. That is, 
$$
\widehat{\boldsymbol{\theta}} := \arg\max_{\boldsymbol{\theta}} \sum^N_{i=1} \log\{\widehat{f}_{T|s_{\boldsymbol{Z},\boldsymbol{\theta}}}(T_i|s_{\boldsymbol{Z},\boldsymbol{\theta}}(\boldsymbol{Z}_i)\}\,.
$$
For example, we can take the kernel density estimator 
$$
\widehat{f}_{T|s_{\boldsymbol{Z},\boldsymbol{\theta}}}(t|s) = \frac{\sum^{N}_{i=1}\mathcal{K}_{h_{1}}(T_i-t)\mathcal{K}_{h_{2}}\{s_{\boldsymbol{Z},\boldsymbol{\theta}}(\boldsymbol{Z}_i)-s\}}{\sum^{N}_{i=1}\mathcal{K}_{h_{2}}\{s_{\boldsymbol{Z},\boldsymbol{\theta}}(\boldsymbol{Z}_i)-s\}}\,,
$$
where $\mathcal{K}_h(x)=\mathcal{K}(x/h)h$ with $K$ a pre-specified kernel function, and $h_1, h_2$ are the rule of thumb bandwidths for $T$ and $s_{\boldsymbol{Z},\boldsymbol{\theta}}(\boldsymbol{Z})$, respectively.
Then with $\widehat{\boldsymbol{\theta}}$ and the moment equation~\eqref{eq:idIndex}, we can estimate $\pi_{\boldsymbol{Z}}$ using the same way as in section~\ref{sec:point estimation}. Since the index model $s_{\boldsymbol{Z},\boldsymbol{\theta}}$ maps $\boldsymbol{Z}$ to a univariate space, both the nonparametric estimator of $f_{T|s_{\boldsymbol{Z},\boldsymbol{\theta}}}(t|s)$ and the sieve approximation of \eqref{eq:idIndex} are two-dimensional nonparametric estimators, reducing the curse of dimensionality.
\end{rk}
\section{Variance Estimation}\label{sec:variance}
Using the expression \eqref{Vtt'}, we propose a plug-in estimator for $V_{tt'}$. For $i=1,\ldots,N$ and $\delta \geq 0$, let $\widehat{d}_{K_0,i}$ defined in section~\ref{sec:B1} in the Appendix be an estimator of $d_{K_0,i}$ defined in \eqref{def:dki}.
Then, a consistent estimator of $V_{tt'}$ is given by
$$\widehat{V}_{tt'}:=u_{K_0}^{\top}(t) \widehat{\Phi}_{K_0\times K_0}^{-1}\cdot \left[\frac{1}{N}\sum_{i=1}^{N}\widehat{d}_{K_0,i}(T,\boldsymbol{M},\boldsymbol{X},Y;\delta)\widehat{d}^{\top}_{K_0,i}(T,\boldsymbol{M},\boldsymbol{X},Y;\delta)\right]\cdot \widehat{\Phi}_{K_0\times K_0}^{-1} u_{K_0}(t), $$
where $\widehat{\Phi}_{K_0\times K_0}=N^{-1}\sum_{i=1}^{N}u_{K_0}(T_i)u^{\top}_{K_0}(T_i).$
From Proposition \ref{rate_pi} and Theorem \ref{thm:sieve}, we have $\sup_{(t,\boldsymbol{z})\in\mathcal{T}\times\mathcal{Z}}|\widehat{\pi}_{K_{\boldsymbol{Z}}}(t,\boldsymbol{z})-\pi_{\boldsymbol{Z}}(t,\boldsymbol{z}) |=o_p(1)$ and $ |\widehat{\mu}(t,t)-\mu(t,t) |\to 0.$ With these results, the consistency of $\widehat{V}_{tt'}$ follows from standard arguments in \cite{chen2007large}.

\section{Selecting the Smoothing Parameters}\label{sec:tunning}
The proposed estimator $\widehat{\mu}(t,t')$ (resp. $\widehat{\mu}_h(t,t')$) in \eqref{def:CBS} (resp. \eqref{def:CBK}) involves tuning parameters $k_1, k_{\boldsymbol{X}},k_{\boldsymbol{M,X}}$, and $K_0$ (resp. $h$). If the treatment $T$ is discrete, $k_1$ and $K_0$ can be determined in the way described in Remark~1, and the kernel type estimator $\widehat{\mu}_h(t,t')$ is not applicable. When the confounders $\boldsymbol{X}$ or the mediators $\boldsymbol{M}$ are discrete, the corresponding sieve bases and parameters $k_{\boldsymbol{X}}$ and $k_{\boldsymbol{M,X}}$ can also be determined in that way. Thus, in this section, we focus on proposing a data-driven method for choosing the parameters for continuous $T, \boldsymbol{X}$ and $\boldsymbol{M}$.

Note that our estimators $\widehat{\mu}(t,t')$ and $\widehat{\mu}_h(t,t')$ for continuous $T$ are nonparametric regression type estimators. For such estimators, the smoothing parameters are usually selected by minimizing certain cross-validation (CV) criteria that approximates the mean squared error (MSE) of the estimator. Although simultaneously selecting all the smoothing parameters using one CV criteria can better approximate the MSE, it is too time-consuming, given so many parameters. Thus, we propose to choose them separately.

First, we propose a method to choose $K_{\boldsymbol{Z}}$ for $\widehat{\pi}_{K_{\boldsymbol{Z}}}$, where $\boldsymbol{Z}\in\{\boldsymbol{X},(\boldsymbol{M},\boldsymbol{Z})\}$, inspired by the least square CV idea for choosing the bandwidth of a kernel density estimator (see e.g., \citealp{fan1996local}). Notice that a weighted integrated squared error (WISE) of the estimated weighting functions can be written in the following way:
\begin{align*}
\int \{\widehat{\pi}_{K_{\boldsymbol{Z}}}(T,\boldsymbol{Z}) - \pi_{\boldsymbol{Z}}(T,\boldsymbol{Z})\}^2 &f_{T,\boldsymbol{Z}}(t,\boldsymbol{z})\,dt\,d\boldsymbol{z} = \int \{\widehat{\pi}_{K_{\boldsymbol{Z}}}(T,\boldsymbol{Z})\}^2 f_{T,\boldsymbol{Z}}(t,\boldsymbol{z})\,dt\,d\boldsymbol{z} \\
&- 2\int \widehat{\pi}_{K_{\boldsymbol{Z}}}(t,\boldsymbol{z}) \pi_{\boldsymbol{Z}}(t,\boldsymbol{z}) f_{T,\boldsymbol{Z}}(t,\boldsymbol{z})\,dt\,d\boldsymbol{z} + \mathbb{E}\left[ \{\pi_{\boldsymbol{Z}}(T,\boldsymbol{Z})\}^2\right]\,,
\end{align*}
where by the definition of $\pi_{\boldsymbol{Z}}$, we have $ \pi_{\boldsymbol{Z}}(t,\boldsymbol{z}) f_{T,\boldsymbol{Z}}(t,\boldsymbol{z}) = f_{T}(t)f_{\boldsymbol{Z}}(\boldsymbol{z})$. This result suggests that minimizing the WISE w.r.t. the number of sieve basis $k_1$ and $k_{\boldsymbol{Z}}$ is equivalent to minimizing
$$
\int \{\widehat{\pi}_{K_{\boldsymbol{Z}}}(T,\boldsymbol{Z})\}^2 f_{T,\boldsymbol{Z}}(t,\boldsymbol{z})\,dt\,d\boldsymbol{z} - 2\int \widehat{\pi}_{K_{\boldsymbol{Z}}}(t,\boldsymbol{z}) f_{T}(t)f_{\boldsymbol{Z}}(\boldsymbol{z})\,dt\,d\boldsymbol{z}\,.
$$

With a finite sample, we can then choose $k_1$ and $k_{\boldsymbol{Z}}$ by minimizing the following least square CV, which is an empirical fully accessible analog of the above criteria
\begin{equation*}
CV(k_1,k_{\boldsymbol{Z}}) = \frac{1}{N}\sum^N_{i=1} \{\widehat{\pi}_{K_{\boldsymbol{Z}}}(T_i,\boldsymbol{Z}_i)\}^2 - \frac{2}{N(N-1)}\sum^{N}_{i=1}\sum^N_{j\neq i} \widehat{\pi}_{K_{\boldsymbol{Z}}}(T_i,\boldsymbol{Z}_j)\,,
\end{equation*}
over a candidate set of $k_1$ and $k_{\boldsymbol{Z}}$. 


To select the parameter $K_0$ in $\widehat{\mu}(t,t')$, we rewrite $\widehat{\mu}(t,t')$ as $\widehat{\mu}_{K_0}(t,t')$ for now. Then, we obtain $K_0$ by using a leave-one-out CV (see \citet[Section 15.2]{li2007nonparametric}), 
\begin{align*}
\widehat{K}_0
=&\arg\min_{K_0}\left[\frac{1}{N}\sum_{i=1}^N \left\{\frac{\widehat{\pi}_{\widehat{K}_{\boldsymbol{M,X}}}(T_i,\boldsymbol{M}_i,\boldsymbol{X}_i)\widehat{\pi}_{\widehat{K}_{\boldsymbol{X}}}(T_i+\delta,\boldsymbol{X}_i)Y_i}{\widehat{\pi}_{\widehat{K}_{\boldsymbol{M,X}}}(T_i+\delta,\boldsymbol{M}_i,\boldsymbol{X}_i)} -\widehat{\mu}^{(-i)}_{K_0}(T_i,T_i+\delta) \right\}^2\right],
\end{align*}
where $\widehat{K}_{\boldsymbol{Z}}=\widehat{k}_1\cdot \widehat{k}_{\boldsymbol{Z}}$ for $\boldsymbol{Z}\in\{\boldsymbol{X},(\boldsymbol{M,X})\}$, and $\widehat{\mu}^{(-i)}_{K_0}(T_i,T_i+\delta)$ is computed as $\widehat{\mu}_{K_0}(T_i,T_i+\delta)$ but without using $\{T_i,\boldsymbol{X}_i,\boldsymbol{M}_i,Y_i\}$.
For $\widehat{\mu}_h(t,t')$, we need to choose the bandwidth $h$. We can also apply a leave-one-out CV method for choosing $h$:
\begin{align*}
\widehat{h}
=&\arg\min_{h}\left[\frac{1}{N}\sum_{i=1}^N \left\{\frac{\widehat{\pi}_{\widehat{K}_{\boldsymbol{M,X}}}(T_i,\boldsymbol{M}_i,\boldsymbol{X}_i)\widehat{\pi}_{\widehat{K}_{\boldsymbol{X}}}(T_i+\delta,\boldsymbol{X}_i)Y_i}{\widehat{\pi}_{\widehat{K}_{\boldsymbol{M,X}}}(T_i+\delta,\boldsymbol{M}_i,\boldsymbol{X}_i)} -\widehat{\mu}^{(-i)}_{h}(T_i,T_i+\delta) \right\}^2\right],
\end{align*}
where $\widehat{\mu}^{(-i)}_{h}(T_i,T_i+\delta)$ is computed as the series estimator $\widehat{\mu}_{h}(T_i,T_i+\delta)$ but without using $\{T_i,\boldsymbol{X}_i,\boldsymbol{M}_i,Y_i\}$.

Apart from the CV method, the rule of thumb bandwidth for kernel methods \cite[see e.g.][]{1986Density} is also available, which sacrifices a small amount of accuracy but is less time-consuming. We experimented with the rule of thumb bandwidth in our numerical studies and obtained satisfactory results. Specifically, the rule of thumb bandwidth $h=C\cdot\text{sd}(T)\cdot N^{-1/5}$, where $\text{sd}(T)$ is the standard deviation of $T$ and $C$ is a constant depending on the kernel function. For example, $C=2.34$ for second-order Epanechnikov kernels, $C=3.03$ for fourth-order Epanechnikov kernels, and $C=1.06$ for the Gaussian kernel. However, Theorem~\ref{thm:ker} requires an undersmoothing bandwidth such that $Nh^5\rightarrow 0$ as $N\rightarrow\infty$. Thus, following the suggestion of \cite{hsu2018direct}, we take $h=C\cdot N^{-0.25}$. 
\section{Monte Carlo Simulation}\label{sec:simulation}
In this section, we conduct Monte Carlo simulations to evaluate the finite sample performance of our proposed estimators.

\subsection{Continuous Treatment}
We consider data-generating processes (DGPs) similar to
the designs in \cite{hsu2018direct}. The confounder $\boldsymbol{X}$ is drawn from the uniform distribution over $[-1.5, 1.5]$. We consider the continuous treatment variable generated from $T=0.3\boldsymbol{X}+\epsilon$, where $\epsilon$ is drawn from the uniform distribution over $[-2,2]$. 
We further generate $U$ and $V$ from the uniform distribution over $[-2,2]$ such that $\epsilon$, $U$, and $V$ are independent of each other. The mediator variable is then generated according to $\boldsymbol{M}=0.3T+0.3\boldsymbol{X}+V$. The outcomes in each scenario are given as follows:
\begin{itemize}
\item Scenario I: $Y = 0.3T+0.3\boldsymbol{M} + 0.5 T\boldsymbol{M} + 0.3 \boldsymbol{X}+ U$.
\item Scenario II: $Y = 0.3T+0.3\boldsymbol{M} + 0.3 \boldsymbol{X}+0.25 T^3 + U$.
\item Scenario III: $Y = 0.3T+0.3M + 0.5 T\boldsymbol{M} + 0.3 \boldsymbol{X}+0.25 T^3 + U$.
\end{itemize}

In Scenario I, there is an interaction effect between $T$ and $\boldsymbol{M}$, with the outcome model being linear. Scenario II considers a nonlinear outcome model; however, there is no interaction effect, implying that the direct and indirect effects are homogeneous, that is, $\mu(t,t)-\mu(t',t)=\mu(t,t')-\mu(t',t')$ and $\mu(t,t)-\mu(t,t')=\mu(t',t)-\mu(t',t')$. In Scenario III, there is an interaction effect, with the outcome model being nonlinear. The true dose-response functions are
\begin{itemize}
\item Scenario I: $\mu(t,t') = 0.3t+0.09t'+0.15 tt'$.
\item Scenario II: $\mu(t,t') = 0.3t+0.09t'+0.25 t^3$.
\item Scenario III: $\mu(t,t') = 0.3t+0.09t'+0.15 tt'+0.25 t^3$.
\end{itemize} 

In all the scenarios, we set the sample size $N=500$ and $1000$. The Monte Carlo trials are repeated for $500$ times. We set $t'=0$ and let $t$
vary over the interval $[-1.5,0)\cup (0,1.5]$. Specifically, we set $t$ equal to the grid
points $\mathcal{G}_r:=\{-1.5,-1.4,...,-0.1, 0.1,...,1.4,1.5\}$. 

To evaluate the performance, we compare the proposed estimators, the covariate-balancing series (CBS) estimators in \eqref{def:CBS} and the covariate-balancing kernel (CBK) regression estimators in \eqref{def:CBK}, with the alternative estimators established in the literature. Specifically, we estimate the weighting functions in both CBS and CBK estimators using the power series. The smoothing parameters are selected by the data-driven method described in Section \ref{sec:tunning}. 

The alternative competitors considered here are the nonparametric weighting kernel (NWK) and semiparametric weighting kernel (SWK) estimators produced by \cite{hsu2018direct} and the linear ordinary least squares (OLS) regression estimators used in the simulation study of \citet[][Section~5]{hsu2018direct}. In particular, SWK (incorrectly) assumes the conditional distributions of $T$ given $\boldsymbol{X}$ and $T$ given $(\boldsymbol{M},\boldsymbol{X})$ are both normal and estimates the distribution parameters using maximum likelihood. OLS estimates the direct and indirect treatment effects by linearly regressing the observed mediator on the observed treatment and covariates and linearly regressing the observed outcome on the observed mediator, treatment and covariates, respectively. Thus, for nonlinear models, OLS is biased.
\begin{table}[tph]
\caption{$10^3\times$ ARMSE for the estimated direct and indirect effects under Scenarios~I -- III, with the smallest value and the smallest nonparametric value in each configuration highlighted in boldface and underline, respectively.}
\label{table:simu}
\par
\begin{center}
\resizebox{\textwidth}{!}{ \begin{tabular}{c|ccc|ccc|ccc|crc}
\hline\hline 
& \multicolumn{6}{c|}{Average Natural Direct Effects }&\multicolumn{6}{c}{Average Natural Indirect Effects}\\
& \multicolumn{3}{l|}{$\text{ARMSE}\{\widehat{\mu}(t,t)-\widehat{\mu}(0,t)\}$}& \multicolumn{3}{l|}{$\text{ARMSE}\{\widehat{\mu}(t,0)-\widehat{\mu}(0,0)\}$} &\multicolumn{3}{l|}{$\text{ARMSE}\{\widehat{\mu}(t,t)-\widehat{\mu}(t,0)\}$} &\multicolumn{3}{l}{$\text{ARMSE}\{\widehat{\mu}(0,t)-\widehat{\mu}(0,0)\}$}\\
&I&II&III&I&II&III&I&II&III&I&II&III\\ \hline
&\multicolumn{12}{c}{N=500}\\ \hline
OLS & 133.70 & 283.79 & 305.15 & {\bf 46.14} & 283.79 & 285.59 & 125.65 & {\bf 15.96} & 125.69 & {\bf 18.17} & {\bf 15.96} & {\bf 18.39} \\ 
SWK &100.79 &122.31 &130.11 &95.14 &123.63 &127.52 &32.85 &31.17 &39.48 & 22.68 &24.25 &24.82 \\
NWK &73.31 &149.91 &155.49 &97.92 &153.31 &169.17 &118.38 &102.31 &119.75 &35.79 &34.79 &34.71 \\ 
CBS &\bf 64.08 &\bf{61.16} &\bf {{70.24}} &57.16 &\bf {54.25} &\bf 62.44 &\bf{32.09} &27.00 &\bf 33.99 &28.02 &27.81 &28.65\\ 
CBK &128.29 &119.89 &130.14 &125.88 &120.57&127.94 &34.98&28.00  &35.78  & 27.01 &26.52 &27.06  \\ \hline
&\multicolumn{12}{c}{N=1000}\\ \hline
OLS &129.72 &281.73 &302.32 &{\bf 34.78} &281.73 &283.04 &124.79 &{\bf 11.05} &124.84 & {\bf 12.48}& {\bf 11.05} & {\bf 12.85} \\ 
SWK &86.19 &87.84 &94.66 &84.07 &90.01 &93.54 &23.54 &22.32 &28.04 &15.38 &15.59 &16.02 \\ 
NWK &59.97 &116.40 &120.46 &90.29 &121.97 &139.49 &107.25 &91.49 &107.44 &29.74 &29.05 &29.12 \\ 
CBS &\bf{45.72} &\bf {41.48} &\bf{49.49} &40.53 &\bf{37.62} &\bf{44.34} &\bf{21.75} &18.44 &\bf{23.22} & 19.18  &18.76 &19.70  \\ 
CBK & 107.37 & 100.23 & 108.27 & 105.57 &100.14 & 106.27 &24.97  & 20.75 & 25.34 &20.19  &20.22  &20.21 \\ \hline
\end{tabular}
} 
\end{center}
\end{table}

For each $t\in \mathcal{G}_r$, we measure the performance of the estimators by the square root mean square errors (Rmse) over the 500 Monte Carlo trials.
Tables~\ref{table:simu} reports the averages of Rmse (ARMSE) for the estimated direct and indirect effects under Scenarios~I to III, where the averages are taken across all the treatment values $t\in \mathcal{G}_r$. Specifically, for any estimator $\widehat{\mu}(t)$ of a function $\mu(t)$ for $t\in \mathcal{G}_r$, $\text{ARMSE}\{\widehat{\mu}(t)\} = |\mathcal{G}_r|^{-1}\sum_{t\in \mathcal{G}_r}\sqrt{ \sum^{500}_{j=1} \{\widehat{\mu}_j(t)-\mu(t)\}^2/500}$, where $|\mathcal{G}_r|$ is the number of elements in the set $\mathcal{G}_r$ and $\widehat{\mu}_j(t)$ denotes the estimate calculated from the $j$th Monte Carlo trial.

Overall, our CBS estimator consistently outperforms the nonparametric NWK across all scenarios, and our proposed CBK method surpasses NWK in every instance except one -- the average natural direct effects under Scenario I. These observations align well with our theories and validate the robust finite sample properties of our methods.Notably, CBS tends to perform better than CBK in most of the cases. The advantage is more obvious for the average natural direct effects, where the models involves non-linear components. This could be attributed to the fact that the tuning parameters for CBS are data-driven, whereas the bandwidth for the kernel function in CBK is determined by a rule-of-thumb approach. 

Note that the results of SWK and NWK are consistent with those in \cite{hsu2018direct}, where the misspecified SWK is better than NWK in most cases. This maybe because the conditional normal distribution assumption can capture the linear relationships of $T$ with $\boldsymbol{X}$ and $(\boldsymbol{M},\boldsymbol{X})$ in the simulated models fairly well. However, our proposed CBS performs better than SWK in most of the cases. 

In particular, under Scenario I, the direct and indirect effects are heterogeneous. Specifically, $\mu(t,0)-\mu(0,0)$ and $\mu(0,t)-\mu(0,0)$ are reduced to linear models, whereas $\mu(t,t)-\mu(0,t)$ and $\mu(t,t)-\mu(t,0)$ are non-linear. It is thus not surprising that the OLS estimator gives the best estimation for $\mu(t,0)-\mu(0,0)$ and $\mu(0,t)-\mu(0,0)$, but performs worst for the others two. 

In Scenario~II, the direct and indirect effects are homogeneous, and the direct effects models are non-linear. The OLS estimators for the direct effects thus perform the worst and do not show any consistency. However, they estimate the indirect effects consistently, and the best since the true indirect effects are linear. Both the proposed CBS and CBK estimators produce smaller ARMSE than NWK. While CBK is comparable to SWK, CBS outperforms SWK.

Under Scenario~III, our proposed CBS and CBK outperform the nonparametric NWK for all the effects. Moreover, the direct effects and indirect effect $\mu(t,t)-\mu(t,0)$ are non-linear and the indirect effect $\mu(0,t)-\mu(0,0)$ is linear. The OLS estimators are not consistent except for $\mu(0,t)-\mu(0,0)$. In those non-linear configurations, our proposed CBS estimator provides the best performance among all the estimators. 
\begin{table}[tph]
\caption{$10^3\times$ ARMSE for the estimated direct and indirect effects under the binary treatment model, with the smallest value in each configuration highlighted in boldface.}
\label{table:simu4}
\par
\begin{center}
\resizebox{\textwidth}{!}{ 
\begin{tabular}{c|c|c|c|c}
\hline\hline
& \multicolumn{2}{c|}{Average Natural Direct Effects }&\multicolumn{2}{c}{Average Natural Indirect Effects}\\
& $\text{ARMSE}\{\widehat{\mu}(1,0)-\widehat{\mu}(0,0)\}$& $\text{ARMSE}\{\widehat{\mu}(1,1)-\widehat{\mu}(0,1)\}$ &$\text{ARMSE}\{\widehat{\mu}(1,1)-\widehat{\mu}(1,0)\}$ &$\text{ARMSE}\{\widehat{\mu}(0,1)-\widehat{\mu}(0,0)\}$\\
\hline
&\multicolumn{4}{c}{n=500} \\ \hline
IPW & 123.59 & 125.53& {94.83} & 43.72 \\ 
CBS & \bf 122.58 &   \bf 122.95 &  \bf 91.53 &  \bf 38.84 \\ \hline
&\multicolumn{4}{c}{n=1000}\\ \hline
IPW & 85.95 & 85.16 & 66.80 & 30.22\\ 
CBS &  \bf 81.03 & \bf 80.71 &  \bf 66.66 &  \bf 28.88\\ \hline 
\end{tabular}
} 
\end{center}
\end{table}
\subsection{Binary Treatment} 
We also conduct a simulation to investigate the finite sample performance when the treatment is binary, that is, $T\in\{0,1\}$. Following \cite{hsu2018nonparametric}, we generate $T$ from a conditional Bernoulli distribution with a success probability of $\exp(\boldsymbol{X})/\{1+\exp(\boldsymbol{X})\}$. The confounder $\boldsymbol{X}$ is generated from the uniform distribution over $[-1.5,1.5]$, the mediator variable is generated through $\boldsymbol{M}=0.3T+0.3\boldsymbol{X}+V$, and the outcome is generated through $Y = 0.3T+0.3\boldsymbol{M} + 0.5 T\boldsymbol{M} + 0.3 \boldsymbol{X}+0.25 T^3 + U$, where $U$ and $V$ are independently drawn from the uniform distribution over $[-2,2]$. The population means of the potential outcomes are $\mu(0,0)=0$, $\mu(0,1)=0.09$, $u(1,0)=0.55$ and $\mu(1,1)=0.79$. 

We compare our proposed CBS estimator with the inverse probability weighting estimator (IPW) introduced by \cite{hsu2018nonparametric}, 
where the generalized propensity score functions $\mathbb{P}(T=t|\boldsymbol{M},\boldsymbol{X})$ and $\mathbb{P}(T=t|\boldsymbol{X})$ are estimated using the nonparametric series logit regression. Table~\ref{table:simu4} reports the ARMSE for the estimated direct and indirect effects. The results show that both estimators improve with the sample size increases, whereas our CBS estimator performs slightly better in all configurations.

\begin{figure}[h]
    \centering
    \par
    {\ \includegraphics[width = .35\textwidth]{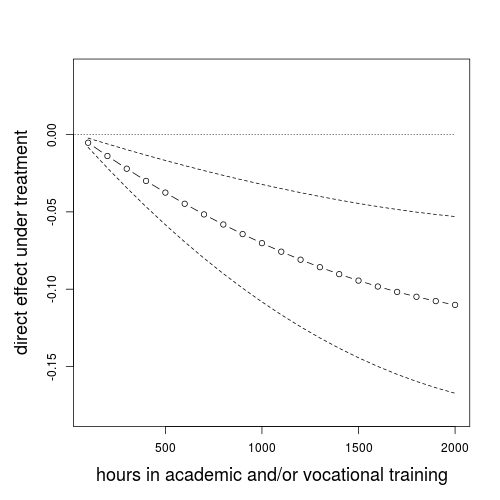}
    \includegraphics[width = .35\textwidth]{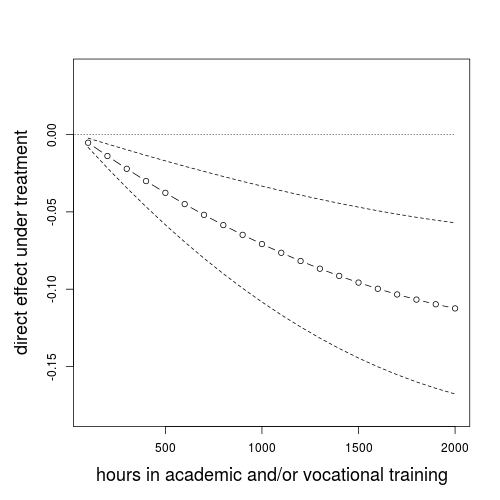}} 
    \par
    {\ \includegraphics[width = .35\textwidth]{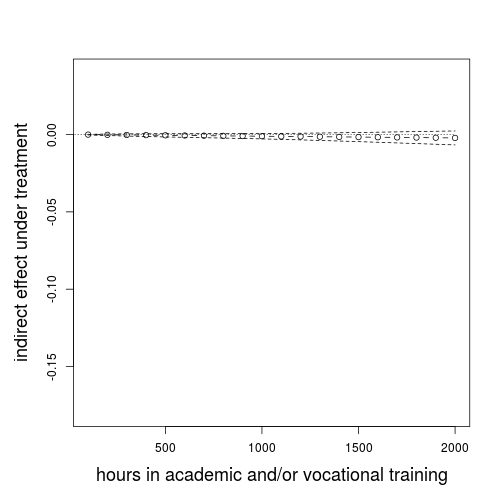}
    \includegraphics[width = .35\textwidth]{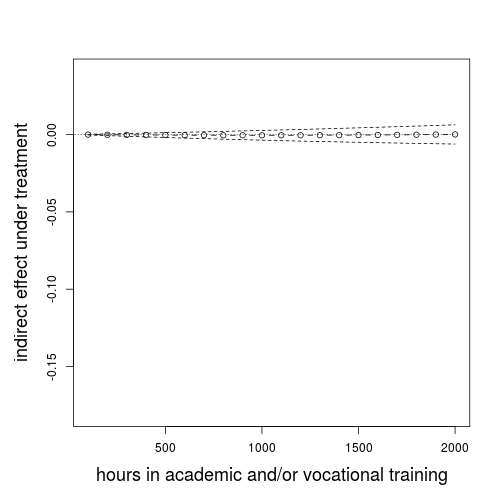}} 
    \par
    \caption{Estimated direct effects $\widehat{\mu}(t,t)-\widehat{\mu}(40,t)$ (top left) and $\widehat{\mu}(t,40)-\widehat{\mu}(40,40)$ (top right), and indirect effects $\widehat{\mu}(t,t)-\widehat{\mu}(t,40)$ (bottom left) and $\widehat{\mu}(40,t)-\widehat{\mu}(40,40)$ (bottom right) for $t\in\{100,200,\cdots,2000\}$, with the estimated $95\%$ confidence bands (dashed lines).} \label{fig:curveplot1}
    \end{figure}

\section{Application} \label{sec:application}
To evaluate the practical value of our method, we revisit the case study of Job Corps analyzed by \cite{schochet2008does} and \cite{hsu2018direct}. Job Corps is a publicly funded U.S. training program that targets economically disadvantaged youths between the ages of 16 and 24 who are legal US residents. Participants received approximately 1200 hours of vocational training and education, housing, and boarding over an average duration of 8 months. Previous analyses by \cite{huber2014identifying}  and \cite{frolich2017direct}  focused on the effect of the program's participation, as a binary treatment effect, on the health and earnings, respectively.
\cite{schochet2008does} find that participation in the Job Corp program increases educational attainment and reduces criminal activity.  
\cite{hsu2018direct} apply a continuous treatment effect analysis to investigate how the total hours spent in either academic or vocational classes during the 12 months affects the participants' criminal activities, mediated through the employment status after the training.

\begin{figure}[h]
    \centering
    \par
    {\ \includegraphics[width = .35\textwidth]{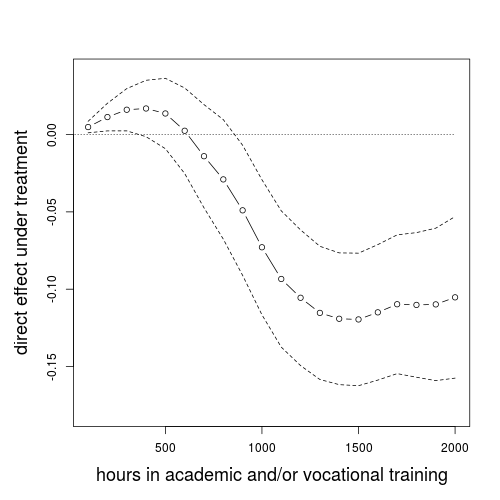}
    \includegraphics[width = .35\textwidth]{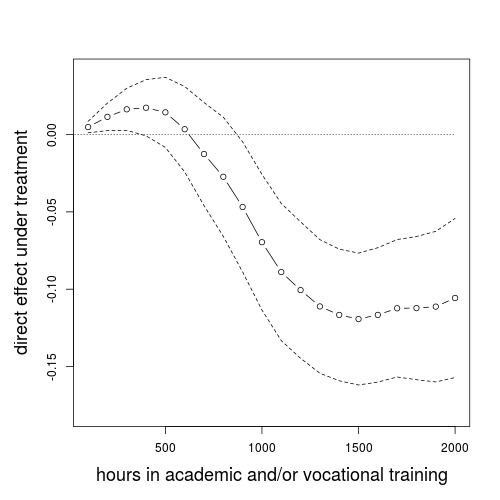}} 
    \par
    {\ \includegraphics[width = .35\textwidth]{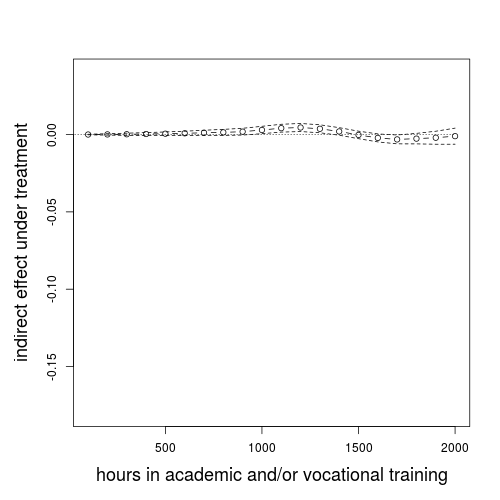}
    \includegraphics[width = .35\textwidth]{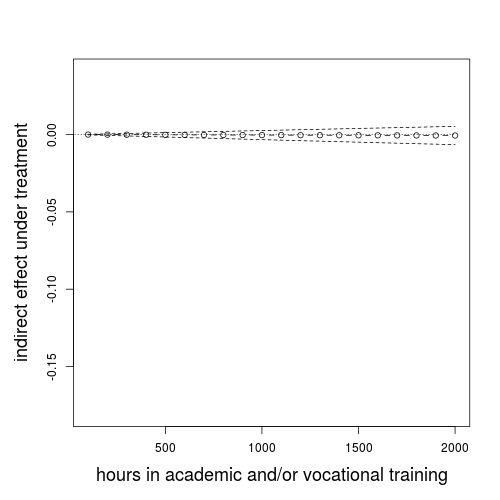}} 
    \par
    \caption{Estimated direct effects $\widehat{\mu}_h(t,t)-\widehat{\mu}_h(40,t)$ (top left) and $\widehat{\mu}_h(t,40)-\widehat{\mu}_h(40,40)$ (top right), and indirect effects $\widehat{\mu}_h(t,t)-\widehat{\mu}_h(t,40)$ (bottom left) and $\widehat{\mu}_h(40,t)-\widehat{\mu}_h(40,40)$ (bottom right) for $t\in\{100,200,\cdots,2000\}$, with the estimated $95\%$ confidence bands (dashed lines).}\label{fig:curveplot2}
\end{figure}

Similar to \cite{hsu2018direct}, our research explores the training program's effects on criminal activity, mediated by post-program employment. But different from them, we first apply our method to assess the direct and indirect effects of attending the program (binary treatment). Subsequently, we utilize our method to examine how these effects vary with the length of participation (continuous treatment) in the program.
Specifically, the mediator variable, denoted by $\boldsymbol{M}$, is the proportion of weeks employed in the second year after the training. The outcome observation, denoted by $Y$, is the number of times an individual was arrested by the police in the fourth year after the program. The confounding variables, denoted by $\boldsymbol{X}$, include a rich set of pre-treatment covariates: age, gender, ethnicity, language competency, education, marital status, household size, and income, previous receipt of social aid, family background (e.g., parents' education), health and health-related behavior at baseline, the expectations about the Job Corps program, and the interaction with the recruiters. 

Let $T$ denote the total hours a participant spent in either academic or vocational classes during the 12-months training program. To analyze the binary treatment mediation effect, we set the binary treatment $D=1$ for all $T>0$ and $D = 0$ for $T=0$. The dataset contains 10,775 individuals, whose post-treatment variables $\boldsymbol{M}$ and $Y$ were fully observed in the follow-up surveys after 2 and 4 years, respectively (see \citealp{hsu2018direct} for a detailed description of the dataset and related statistics). The empirical results, as summarized in Table \ref{table:emprical}, indicate that both direct and indirect effects of this treatment on reducing criminal behavior are not significant. This lack of significance suggests that simply attending the training, regardless of the duration, does not have a causal effect on altering criminal behavior. 
Given the finding, it is important to consider not just the presence of the treatment (i.e., attending the training) but also its intensity (how many hours were spent in training).

\begin{table}
\caption{$10\times$ Direct and indirect effects with confidence interval under the binary treatment}
\label{table:emprical}
\par
\begin{center}
	\resizebox{\textwidth}{!}{ 
		\begin{tabular}{c|c|c|c|c}
			\hline\hline
			& \multicolumn{2}{c|}{Average Natural Direct Effects }&\multicolumn{2}{c}{Average Natural Indirect Effects}\\
			& $\{\widehat{\mu}(1,0)-\widehat{\mu}(0,0)\}$& $\{\widehat{\mu}(1,1)-\widehat{\mu}(0,1)\}$ &$\{\widehat{\mu}(1,1)-\widehat{\mu}(1,0)\}$ &$\{\widehat{\mu}(0,1)-\widehat{\mu}(0,0)\}$\\
			\hline
		
			CBS & -0.075 (-0.348, 0.184) & -0.082 (-0.336,0.185)& 0.005 (-0.041,0.051) &0.011 (-0.050,0.072) \\ 
		 \hline 
		\end{tabular}
	} 
\end{center}
\end{table}

We then analyze the 4,000 individuals, from the 10,775 ones in the original sample, who received a positive treatment intensity, that is, $T>0$. We apply the proposed method to this dataset to analyze the mediation effects to help policymakers design more efficient intervention programs. The benchmark treatment level is fixed at $t'=40$, that is, a rather small intensity of $40$ hours. We compute the response curve estimators $\widehat{\mu}(t,t')$ and $\widehat{\mu}_h(t,t')$ proposed in \eqref{def:CBS} and \eqref{def:CBK}, respectively. Then, we use them to obtain the estimates of the direct and indirect effects for varying $t\in\{100, 200, ..., 1900, 2000\}$. The smoothing parameters are determined using the data-driven approach described in Section \ref{sec:tunning}.

Figure \ref{fig:curveplot1} reports the estimated direct and indirect effects based on our CBS $\widehat{\mu}(t,t')$, while Figure \ref{fig:curveplot2} reports the estimated direct and indirect effects based on our CBK $\widehat{\mu}_h(t,t')$. Figure \ref{fig:curveplot1} shows that the direct effects are significantly negative at the $5\%$ level for all $t$ values considered, with the effects becoming larger as the training time increases. Additionally, all estimated indirect effects are insignificant. 

The results in Figure \ref{fig:curveplot2} are similar to those of \cite{hsu2018direct}. They show that small treatment intensities do not reduce the number of arrests. From $800$ hours onward, the direct effects become significantly negative, the $95\%$ confidence bands exclude zero, and the effect peaks occur around $1400$ hours. The estimated indirect effects are mostly insignificant, consistent with the results presented in Figure \ref{fig:curveplot1}. In conclusion, the Job Corps program has significant direct effects on the number of arrests in the fourth year. 
In contrast, for the investigated range of treatment intensities, the indirect effects of program-induced employment changes on arrests are close to zero. 

Comparing Figure 1 to Figure 2 and the results in \cite{hsu2018direct}, our sieve regression estimator gives smoother results than the kernel-type ones. The kernel-type estimators seem to have some boundary effects that show unstable results on the boundaries of the support of the treatment values. Recall that in our simulation, with or without interaction for both linear and nonlinear models, our CBS gives better direct effects estimation than CBK. Thus, Figure~1 may provide a more reliable result. However, to confirm these analyses, some specification and monotonicity tests are needed, which can be an interesting future direction.

\section{Discussion and Conclusion \label{sec:conclusion}}
This study provides a novel approach for estimating causal mediation effects, which unifies the binary, multi-valued, and continuous treatments and the mixture of discrete and continuous treatments under a sequential ignorability condition. Furthermore, we establish the asymptotic normality for our proposed estimators. In particular, we show that our proposed estimators attain the semiparametric efficiency bounds when the treatment is discrete and asymptotically more efficient than the existing method when the treatment is continuous.
Nevertheless, our study has some limitations that should be addressed by future research. First, an extension to allow for high-dimensional covariates is needed since the fully nonparametric estimation suffers from the curse of dimensionality. Second, our idea is readily applicable to panel data, which is commonly observed in practice; however, the asymptotic analysis could be more difficult and a worthwhile future study.

\section*{Acknowledgments}
The authors sincerely thank the editor Esfandiar Maasoumi and the referees for their
constructive suggestions and comments. Wei
Huang's research is supported by the Professor Maurice H. Belz Fund of the
University of Melbourne.  Zheng Zhang is supported by the fund from the fund from the National Natural Science Foundation of Beijing,
China [grant number 1222007] and the fund for building world-class universities (disciplines) of Renmin University of China
[project number KYGJC2023011].

\bigskip
\begin{center}
{\large\bf SUPPLEMENTARY MATERIAL}
\end{center}
{\bf Supplementary Material for “Nonparametric Estimation of Mediation Effects with A General Treatment”:} The supplementary material is only for online publication (pdf file). It contains the assumptions required to derive the asymptotic properties of $\widehat{\pi}_{\boldsymbol{Z}}$ and detailed discussions on the assumptions, the asymptotic results of $\widehat{\pi}_{\boldsymbol{Z}}$ and the proofs of Theorems~\ref{thm:sieve}, \ref{thm:mutilde} and Corollary~\ref{cor:efficiency}.

\if0\blind
\if1\blind
{} \fi
\bibliographystyle{agsm}
\bibliography{Semiparametric.bib}
\clearpage\appendix
\section*{Appendix \label{sec:appendix}}
\setcounter{section}{0} 
\section{Proof of \eqref{id:mutt'}}\label{app:id}
We rewrite the identification of $\mu(t,t')$ in \cite{hsu2018direct} as follows:
\begin{align*}
\mu(t,t')=&\mathbb{E}\left[ \frac{f_T(t)}{f_{T|\boldsymbol{M},\boldsymbol{X}}(t|\boldsymbol{M},\boldsymbol{X})} \cdot \frac{f_{T|\boldsymbol{M},\boldsymbol{X}}(t'|\boldsymbol{M},\boldsymbol{X})}{f_{T|\boldsymbol{X}}(t'|\boldsymbol{X})} Y\bigg|T=t\right]\\
=&\mathbb{E}\left[ \frac{f_T(t)}{f_{T|\boldsymbol{M},\boldsymbol{X}}(t|\boldsymbol{M},\boldsymbol{X})} \cdot \frac{f_{T|\boldsymbol{M},\boldsymbol{X}}(t'|\boldsymbol{M},\boldsymbol{X})}{f_{T}(t')}\cdot \frac{f_{T}(t')}{f_{T|\boldsymbol{X}}(t'|\boldsymbol{X})} Y\bigg|T=t\right]\\
=&\mathbb{E}\left[ \frac{f_T(T)}{f_{T|\boldsymbol{M},\boldsymbol{X}}(T|\boldsymbol{M},\boldsymbol{X})} \cdot \frac{f_{T|\boldsymbol{M},\boldsymbol{X}}(T+\delta|\boldsymbol{M},\boldsymbol{X})}{f_{T}(T+\delta)}\cdot \frac{f_{T}(T+\delta)}{f_{T|\boldsymbol{X}}(T+\delta|\boldsymbol{X})} Y\bigg|T=t\right]\\
=&\mathbb{E}\left[\frac{{\color{black}\pi_{\boldsymbol{M},\boldsymbol{X}}}(T,\boldsymbol{M},\boldsymbol{X})}{{\color{black}\pi_{\boldsymbol{M},\boldsymbol{X}}}(T+\delta,\boldsymbol{M},\boldsymbol{X})} \cdot {\color{black}\pi_{\boldsymbol{X}}}(T+\delta,\boldsymbol{X}) Y\bigg|T=t\right].
\end{align*}
\section{Some Preliminary Results}\label{appendix:preliminary}
We recall some preliminary results on the convergence rates of $\widehat{\pi}_{K_{\boldsymbol{Z}}}(t,\boldsymbol{Z})$ that are directly implied by the results in \cite{Ai_Linton_Motegi_Zhang_cts_treat}. 
We impose the following conditions based on those of \cite{Ai_Linton_Motegi_Zhang_cts_treat}: {For $\boldsymbol{Z} \in \{\boldsymbol{X},(\boldsymbol{M},\boldsymbol{X})\}$, we assume
\begin{assumption}
\label{as:suppX} (i) The support of $\boldsymbol{Z}$, $\mathcal{Z}$ is a compact set. (ii) There exist two positive constants $\eta_{1}$ and $\eta_{2}$ such that
\begin{align*}
&0<\eta_{1}\leq \pi_{\boldsymbol{Z}}(t,\boldsymbol{z})\leq\eta_{2}<\infty\ ,\ \forall
(t,\boldsymbol{z})\in\mathcal{T}\times\mathcal{Z}.
\end{align*} 
\end{assumption}
\begin{assumption}
\label{as:smooth_pi} There exist $\Lambda_{k_{1}\times k_{\boldsymbol{Z}}}\in\mathbb{R}^{k_{1}\times k_{\boldsymbol{Z}}}$ and constant $\alpha_{\boldsymbol{Z}}>0$ such that
\[
\sup_{(t,\boldsymbol{z})\in\mathcal{T}\times\mathcal{Z}}\left\vert
(\rho^{\prime-1}\left\{\pi_{\boldsymbol{Z}}(t,\boldsymbol{z})\right\}-u_{k_{1}}(t)^{\top}\Lambda_{k_{1}\times k_{\boldsymbol{Z}}}v_{k_{\boldsymbol{Z}}}(\boldsymbol{z})\right\vert
=O(K_{\boldsymbol{Z}}^{-\alpha_{\boldsymbol{Z}}})\,,
\]
where $\rho^{\prime-1}(v)=-\log v-1$.
\end{assumption}
\begin{assumption}
\label{as:u&v}(i) The eigenvalues of $\mathbb{E}[u_{k_{1}}(T)u_{k_{1}}(T)^{\top}]$, $\mathbb{E}[v_{k_{\boldsymbol{Z}}}(
\boldsymbol{Z})$ $v_{k_{\boldsymbol{Z}}}(\boldsymbol{Z})^{\top}] $ are bounded away from zero and infinity uniformly in $k_{1}$, $k_{\boldsymbol{Z}}$. (ii) There are sequences of
constants $\zeta_{1}(k_{1})$ and $\zeta_{\boldsymbol{Z}}(k_{\boldsymbol{Z}})$ satisfying $\sup_{t\in\mathcal{T}}\Vert u_{k_{1}}(t)\Vert\leq\zeta_{1}(k_{1})$ and $\sup_{\boldsymbol{z}\in\mathcal{Z}}\Vert v_{k_{\boldsymbol{Z}}}(\boldsymbol{z})\Vert\leq\zeta_{\boldsymbol{Z}}(k_{\boldsymbol{Z}})$ such that $\sqrt{N}K_{\boldsymbol{Z}}^{-\alpha_{\boldsymbol{Z}}}\rightarrow 0$ and $\zeta(K_{\boldsymbol{Z}})\sqrt{K_{\boldsymbol{Z}}^{2}/N}\rightarrow 0$ as $N\rightarrow\infty$, where $K_{\boldsymbol{Z}}=k_1k_{\boldsymbol{Z}}$ and $\zeta(K_{\boldsymbol{Z}})=\zeta_{1}(k_{1})\zeta_{\boldsymbol{Z}}(k_{\boldsymbol{Z}})$.
\end{assumption}
}
Assumption \ref{as:suppX} (i) requires the covariates, the treatment
variable and the mediator to be bounded. This condition, despite being restrictive, is commonly
imposed in the non-parametric regression literature. However, we can replace it
with a restriction on the tail distribution of $(\boldsymbol{M},\boldsymbol{X},T)$. For
example, \citet*[Assumption 3]{chen2008semiparametric} assume that the
support of $\boldsymbol{X}$ is the entire Euclidean space, but impose
$\int_{\mathbb{R}^{r}}(1+|\boldsymbol{x}|^{2})^{\omega}f_{\boldsymbol{x}}(\boldsymbol{x})d\boldsymbol{x}<\infty$ for some
$\omega>0$. Assumption \ref{as:suppX} (ii) requires the weighting function to
be bounded and bounded away from zero. We can relax Assumption \ref{as:suppX}
(ii) by allowing $\eta_{1}$ ($\eta_{2}$) to go to zero (infinity) slowly as
$N\rightarrow\infty$. Notice that $u_{k_{1}}(t)^{\top}\Lambda v_{k_{\boldsymbol{z}}%
}(\boldsymbol{z})$ is a linear sieve approximation for $\rho^{\prime-1}\left\{\pi_{\boldsymbol{z}}(t,\boldsymbol{z})\right\}$. Assumption \ref{as:smooth_pi} requires
the sieve approximation error to shrink to zero at a polynomial rate. A
variety of sieve basis functions satisfy this condition. When both $T$ and $\boldsymbol{Z}$ are discrete, $K_{\boldsymbol{Z}}$ is a finite constant independent of $N$ and $\alpha_{\boldsymbol{Z}}=+\infty$. When $(T,\boldsymbol{Z})$ has continuous or mixed components, $K_{\boldsymbol{Z}}\rightarrow\infty$ as $N\rightarrow \infty$ and $\alpha_{\boldsymbol{Z}}$ is positively affected by the smoothness of $\rho^{\prime-1}\left(\cdot \right) $ and negatively affected by the number of continuous and mixed
components. Assumption \ref{as:u&v} (i) ensures that the sieve estimator is
non-degenerate. This condition is common in the sieve regression literature
(see \citealp{andrews1991asymptotic} and \citealp{Newey97}). If the approximation
error is nonzero, Assumption \ref{as:u&v} (ii) imposes a restriction on the
growth rate of the smoothing parameters $k_{1}$ and $k_{\boldsymbol{z}}$ to ensure under-smoothing.
Under these conditions, the following are direct results from \citet[Propositions 1 and 2]{Ai_Linton_Motegi_Zhang_cts_treat}:
\begin{prop}
\label{rate_pi} Suppose that Assumptions \ref{as:suppX} -- \ref{as:u&v} hold. {For $\boldsymbol{Z}\in \{\boldsymbol{X}, (\boldsymbol{M},\boldsymbol{X})\}$, we} have
\[sup_{(t,\boldsymbol{m},\boldsymbol{x})\in\mathcal{T}\times\mathcal{M}\times\mathcal{X} }|\widehat{\pi}_{K_{\boldsymbol{Z}}}(t,\boldsymbol{z})-\pi
_{\boldsymbol{Z}}(t,\boldsymbol{z})|=O_{p}\left(\max\left\{\zeta(K_{\boldsymbol{Z}})K_{\boldsymbol{Z}}^{-\alpha_{\boldsymbol{Z}}},\zeta(K_{\boldsymbol{Z}})\sqrt{\frac{K_{\boldsymbol{Z}}}{N}}\right\}\right)\ ,\]
and 
\[
\int_{\mathcal{T}\times\mathcal{Z}}|\widehat{\pi}_{K_{\boldsymbol{Z}}}(t,\boldsymbol{z})-\pi
_{\boldsymbol{Z}}(t,\boldsymbol{z})|^{2}dF_{T,\boldsymbol{Z}}(t,\boldsymbol{z})=O_{p}\left(
\max\left\{K_{\boldsymbol{Z}}^{-2\alpha_{\boldsymbol{Z}}},\frac{K_{\boldsymbol{Z}}}{N}\right\} \right) \,,
\]
and
\[
\frac{1}{N}\sum_{i=1}^{N}|\widehat{\pi}_{K_{\boldsymbol{Z}}}(T_{i},\boldsymbol{Z}_{i})-\pi_{\boldsymbol{Z}}(T_{i},\boldsymbol{Z}_{i})|^{2}=O_{p}\left( \max\left\{ K_{\boldsymbol{Z}}^{-2\alpha_{\boldsymbol{Z}}},\frac{K_{\boldsymbol{Z}}}{N}\right\} \right)\,.
\]
\end{prop}
\begin{prop}\label{prop:influence_pi} Assume that Assumptions \ref{as:suppX} -- \ref{as:u&v} hold. For for any square-integrable random variable $\phi(T,\boldsymbol{M},\boldsymbol{X},Y)\in L^{2}$ and $\boldsymbol{Z}\in \{\boldsymbol{X}, (\boldsymbol{M},\boldsymbol{X})\}$, if there exist a $\Gamma_{k_1\times k_{\boldsymbol{Z}}}\in \mathbb{R}^{k_1\times k_{\boldsymbol{Z}}}$ and a constant $\gamma_{\boldsymbol{Z}}>0$, s.t. $$\sup_{t \times \boldsymbol{z} \in \mathcal{T} \times \mathcal{Z}}\Big|\mathbb{E}[\phi(T,\boldsymbol{M},\boldsymbol{X},Y)|T=t,\boldsymbol{Z}=\boldsymbol{z}] - u_{k_1}(t)^\top \Gamma_{k_1\times k_{\boldsymbol{Z}}} v_{k_{\boldsymbol{Z}}}(\boldsymbol{z})\Big|=O(K_{\boldsymbol{Z}}^{-\gamma_{\boldsymbol{Z}}})\,,$$ 
then we have, for any $\delta\geq 0$, 
\begin{align}\label{eq:effpi}
& \frac{1}{\sqrt{N}}\sum_{i=1}^{N}\{\widehat{\pi}_{K_{\boldsymbol{Z}}}(T_{i}+\delta,\boldsymbol{Z}_{i})\phi(T_{i},\boldsymbol{M}_i,\boldsymbol{X}_{i},Y_{i})-\mathbb{E}[\pi_{\boldsymbol{Z}}(T+\delta,\boldsymbol{Z})\phi(T,\boldsymbol{M},\boldsymbol{X},Y)]\}\\
=&\frac{1}{\sqrt{N}}\sum_{i=1}^{N} \text{IF}_{\pi_{\boldsymbol{Z}},i}\{\delta,\phi(T,\boldsymbol{M},\boldsymbol{X},Y)\} \notag\\
&\qquad\qquad+O_p\left(\sqrt{N}K_{\boldsymbol{Z}}^{-\alpha_{\boldsymbol{Z}}}\right)+O_p\left(K_{\boldsymbol{Z}}^{-\gamma_{\boldsymbol{Z}}}\right)+O_p\left(\zeta(K_{\boldsymbol{Z}})\sqrt{\frac{K_{\boldsymbol{Z}}^2}{N}} \right)\notag
\end{align}
where, for $i=1,\ldots, N$,
\begin{align}\label{def:IF_piZ}
\text{IF}_{\pi_{\boldsymbol{Z}},i}\{\delta,\phi(T,\boldsymbol{M},\boldsymbol{X},Y)\}=&\pi_{\boldsymbol{Z}}(T_{i}+\delta,\boldsymbol{Z}_{i})\phi(T_{i},\boldsymbol{M}_i,\boldsymbol{X}_{i},Y_{i})\\
&-\pi_{\boldsymbol{Z}}(T_{i},\boldsymbol{Z}_{i})\frac{f_{T|\boldsymbol{Z}}(T_i-\delta|\boldsymbol{Z}_i)}{f_{T|\boldsymbol{Z}}(T_i|\boldsymbol{Z}_i)}\mathbb{E}[\phi(T_{i}-\delta,\boldsymbol{M}_i,\boldsymbol{X}_{i},Y_{i})|T_{i},\boldsymbol{Z}_{i}]\notag\\
&+\mathbb{E}[\pi_{\boldsymbol{Z}}(T_{i},\boldsymbol{Z}_{i})\frac{f_{T|\boldsymbol{Z}}(T_i-\delta|\boldsymbol{Z}_i)}{f_{T|\boldsymbol{Z}}(T_i|\boldsymbol{Z}_i)}\phi(T_{i}-\delta,\boldsymbol{M}_i,\boldsymbol{X}_{i},Y_{i})|\boldsymbol{Z}_{i}]\notag\\
&-\mathbb{E}[\pi_{\boldsymbol{Z}}(T_{i},\boldsymbol{Z}_{i})\frac{f_{T|\boldsymbol{Z}}(T_i-\delta|\boldsymbol{Z}_i)}{f_{T|\boldsymbol{Z}}(T_i|\boldsymbol{Z}_i)}\phi(T_{i}-\delta,\boldsymbol{M}_i,\boldsymbol{X}_{i},Y_{i})]\notag\\
&+\mathbb{E}[\pi_{\boldsymbol{Z}}(T_{i},\boldsymbol{Z}_{i})\frac{f_{T|\boldsymbol{Z}}(T_i-\delta|\boldsymbol{Z}_i)}{f_{T|\boldsymbol{Z}}(T_i|\boldsymbol{Z}_i)}\phi(T_{i}-\delta,\boldsymbol{M}_i,\boldsymbol{X}_{i},Y_{i})|T_{i}]\notag\\
&-\mathbb{E}[\pi_{\boldsymbol{Z}}(T_{i},\boldsymbol{Z}_{i})\frac{f_{T|\boldsymbol{Z}}(T_i-\delta|\boldsymbol{Z}_i)}{f_{T|\boldsymbol{Z}}(T_i|\boldsymbol{Z}_i)}
\phi(T_{i}-\delta,\boldsymbol{M}_i,\boldsymbol{X}_{i},Y_{i})]\,,\notag
\end{align} 
and $\mathbb{E}[\text{IF}_{\pi_{\boldsymbol{Z}},i}\{\delta,\phi(T,\boldsymbol{M},\boldsymbol{X},Y)\}] = 0$.
\end{prop}

Using these results, we show the following lemma that is useful for deriving our main theorems. The proof can be found in section~S2 in the supplemental material.
\begin{lemma}\label{lem:influence}
Assume that Assumptions \ref{as:suppX} -- \ref{as:u&v} hold. For for any square-integrable random variable $\phi(T,\boldsymbol{M},\boldsymbol{X},Y)\in L^{2}$ and $\boldsymbol{Z}\in \{\boldsymbol{X}, (\boldsymbol{M},\boldsymbol{X})\}$, if there exist a $\Gamma_{k_1\times k_{\boldsymbol{Z}}}\in \mathbb{R}^{k_1\times k_{\boldsymbol{Z}}}$ and a constant $\gamma_{\boldsymbol{Z}}>0$, s.t. $$\sup_{t \times \boldsymbol{z} \in \mathcal{T} \times \mathcal{Z}}\Big|\mathbb{E}[\phi(T,\boldsymbol{M},\boldsymbol{X},Y)|T=t,\boldsymbol{Z}=\boldsymbol{z}] - u_{k_1}(t)^\top \Gamma_{k_1\times k_{\boldsymbol{Z}}} v_{k_{\boldsymbol{Z}}}(\boldsymbol{z})\Big|=O(K_{\boldsymbol{Z}}^{-\gamma_{\boldsymbol{Z}}})\,,$$ 
then we have, for any $\delta\geq 0$, 
\begin{align}\label{eq:effrk}
& \frac{1}{\sqrt{N}}\sum_{i=1}^{N}\bigg\{\frac{\widehat{\pi}_{K_{\boldsymbol{M,X}}}(T_i,\boldsymbol{M}_i,\boldsymbol{X}_i)}{\widehat{\pi}_{K_{\boldsymbol{M,X}}}(T_i+\delta,\boldsymbol{M}_i,\boldsymbol{X}_i)}\widehat{\pi}_{K_{\boldsymbol{X}}}(T_i+\delta,\boldsymbol{X}_i)\phi(T_{i},\boldsymbol{M}_i,\boldsymbol{X}_{i},Y_{i})\notag\\
&\qquad \qquad \qquad -\mathbb{E}\left[\frac{\pi_{\boldsymbol{M,X}}(T,\boldsymbol{M},\boldsymbol{X})}{\pi_{\boldsymbol{M,X}}(T_i+\delta,\boldsymbol{M}_i,\boldsymbol{X}_i)}\pi_{\boldsymbol{X}}(T+\delta,\boldsymbol{X})\phi(T,\boldsymbol{M},\boldsymbol{X},Y)\right]\bigg\}\\
=&\frac{1}{\sqrt{N}}\sum_{i=1}^{N}\Bigg[ \text{IF}_{\pi_{\boldsymbol{X}},i}\bigg\{\delta,\frac{\pi_{\boldsymbol{M},\boldsymbol{X}}(T,\boldsymbol{M},\boldsymbol{X})}{\pi_{\boldsymbol{M},\boldsymbol{X}}(T+\delta,\boldsymbol{M},\boldsymbol{X})}\cdot\phi(T,\boldsymbol{M},\boldsymbol{X},Y)\bigg\} \notag\\
&\qquad \qquad +\text{IF}_{\pi_{\boldsymbol{M},\boldsymbol{X}},i}\bigg\{0,\frac{\pi_{\boldsymbol{X}}(T+\delta,\boldsymbol{X})}{\pi_{\boldsymbol{M},\boldsymbol{X}}(T+\delta,\boldsymbol{M},\boldsymbol{X})}\cdot\phi(T,\boldsymbol{M},\boldsymbol{X},Y)\bigg\} \notag\\
&\qquad \qquad -\text{IF}_{\pi_{\boldsymbol{M},\boldsymbol{X}},i}\bigg\{\delta,\frac{\pi_{\boldsymbol{M},\boldsymbol{X}}(T,\boldsymbol{M},\boldsymbol{X})\pi_{\boldsymbol{X}}(T+\delta,\boldsymbol{X})}{\pi^2_{\boldsymbol{M,X}}(T+\delta, \boldsymbol{M,X})}\cdot\phi(T,\boldsymbol{M},\boldsymbol{X},Y)\bigg\} \Bigg] \notag\\
&\qquad+O_p\left(\sqrt{N}K_{\boldsymbol{X}}^{-\alpha_{\boldsymbol{X}}}\right)+O_p\left(K_{\boldsymbol{X}}^{-\gamma_{\boldsymbol{X}}}\right)+O_p\left(\zeta(K_{\boldsymbol{X}})\sqrt{\frac{K_{\boldsymbol{X}}^2}{N}} \right)\notag\\
&\qquad+O_p\left(\sqrt{N}K_{\boldsymbol{M},\boldsymbol{X}}^{-\alpha_{\boldsymbol{M},\boldsymbol{X}}}\right)+O_p\left(K_{\boldsymbol{M},\boldsymbol{X}}^{-\gamma_{\boldsymbol{M},\boldsymbol{X}}}\right)+O_p\left(\zeta(K_{\boldsymbol{M},\boldsymbol{X}})\sqrt{\frac{K_{\boldsymbol{M},\boldsymbol{X}}^2}{N}} \right)\,.\notag
\end{align}
\end{lemma} 

\subsection{Estimating the $d_{K_0,i}$'s}\label{sec:B1}
To estimate the asymptotic variance $V_{tt'}$ in Section~\ref{sec:variance}, for $i=1,\ldots, N$, we estimate $d_{K_0,i}$ by
\begin{align*}
    \widehat{d}_{K_0,i}(T,\boldsymbol{M},\boldsymbol{X},Y;\delta):=&\widehat{\text{IF}}_{\pi_{\boldsymbol{X}},i}\bigg\{\delta,\frac{\widehat{\pi}_{K_{\boldsymbol{M,X}}}(T,\boldsymbol{M},\boldsymbol{X})}{\widehat{\pi}_{K_{\boldsymbol{M,X}}}(T+\delta,\boldsymbol{M},\boldsymbol{X})}\cdot u_{K_0}(T) Y \bigg\} \notag\\
    &+\widehat{\text{IF}}_{\pi_{\boldsymbol{M},\boldsymbol{X}},i}\bigg\{0,\frac{\widehat{\pi}_{K_{\boldsymbol{X}}}(T+\delta,\boldsymbol{X})}{\widehat{\pi}_{K_{\boldsymbol{M,X}}}(T+\delta,\boldsymbol{M},\boldsymbol{X})}\cdot u_{K_0}(T) Y \bigg\} \notag\\
    &-\widehat{\text{IF}}_{\pi_{\boldsymbol{M},\boldsymbol{X}},i}\bigg\{\delta,\widehat{\pi}_{K_{\boldsymbol{M,X}}}(T,\boldsymbol{M},\boldsymbol{X})\widehat{\pi}_{K_{\boldsymbol{X}}}(T+\delta,\boldsymbol{X})\cdot u_{K_0}(T) Y \bigg\}\\
    &-\widehat{\mathbb{E}}\left[\frac{\widehat{\pi}_{K_{\boldsymbol{M},\boldsymbol{X}}}(T_i,\boldsymbol{M}_i,\boldsymbol{X}_i)\widehat{\pi}_{K_{\boldsymbol{X}}}(T_i+\delta,\boldsymbol{X}_i)}{\widehat{\pi}_{K_{\boldsymbol{M,X}}}(T_i+\delta, \boldsymbol{M_i,X_i})}\cdot u_{K_0}(T_i) Y_i|T_i\right]\\
    &+\frac{1}{N}\sum_{i=1}^N\frac{\widehat{\pi}_{K_{\boldsymbol{M},\boldsymbol{X}}}(T_i,\boldsymbol{M}_i,\boldsymbol{X}_i)\widehat{\pi}_{K_{\boldsymbol{X}}}(T_i+\delta,\boldsymbol{X}_i)}{\widehat{\pi}_{K_{\boldsymbol{M,X}}}(T_i+\delta, \boldsymbol{M_i,X_i})}\cdot u_{K_0}(T_i) Y_i\,,
    \end{align*}
    where, for $\boldsymbol{Z}\in \{\boldsymbol{X},(\boldsymbol{M},\boldsymbol{X})\}$, $\widehat{\text{IF}}_{\pi_{\boldsymbol{Z}},i}$ is an estimator of $\text{IF}_{\pi_{\boldsymbol{Z}},i}$ defined in \eqref{def:IF_piZ} in Appendix~\ref{appendix:preliminary}, and $\widehat{\mathbb{E}}$ denotes the least square regression. Specifically, we estimate $\pi_{\boldsymbol{Z}}$, $f_{T|\boldsymbol{Z}}$, the conditional expectations, and expectations in \eqref{def:IF_piZ} by $\widehat{\pi}_{K_{\boldsymbol{Z}}}$, the conditional kernel density estimation, the least square regression of the estimated response variable on the corresponding sieve basis, and the sample average of the estimated variables, respectively. For example, $\mathbb{E}\left[\pi_{\boldsymbol{X}}(T_i,\boldsymbol{X}_i) u_{K_0}(T_i)Y_i|T_i,\boldsymbol{X}_i\right]$ is estimated by the least square regression of $\widehat{\pi}_{K_{\boldsymbol{X}}}(T_i,\boldsymbol{X}_i)u_{K_0}(T_i)Y_i$ on a sieve basis $w_{K_{T,\boldsymbol{X}}}(T_i,\boldsymbol{X}_i)$.

\section{Asymptotics of Kernel Estimators under Continuous Treatments}\label{app:asymptoics_kernel}
To derive the asymptotic normality of our kernel regression estimators, the following assumptions are imposed.
\begin{assumption}
\label{as:K} $\mathcal{K}(\cdot)$ is a univariate kernel function symmetric around the
origin that satisfies (i) $\int \mathcal{K}(u)du=1$; (ii) $\int u^{2}\mathcal{K}(u)du=\kappa
_{21}\in(0,\infty)$; (iii) $\int \mathcal{K}^{2}(u)du=\kappa_{02}<\infty$; and (iv)
$\int|\mathcal{K}(u)|^{2+\delta}du<\infty$, for some $\delta>0$.
\end{assumption}
\begin{assumption}
\label{as:h} As $N\to\infty$, $h\to0$, $Nh\to\infty$ and $Nh^5\to 0$.
\end{assumption}
Assumption \ref{as:h} is common in the
kernel regression literature (see \cite{li2007nonparametric}).
\begin{theorem}
\label{thm:ker} Suppose Assumptions \ref{as:SequentialIgnore}, \ref{as:suppX} -- \ref{as:h} hold. Then, we have
$$\sqrt{Nh}\left\{ \widehat{\mu}_h(t,t')-\mu(t,t')\right\}=\sqrt{\frac{h}{N}}\sum_{i=1}^N\psi_{tt'}(Y_i,T_i,\boldsymbol{M}_i,\boldsymbol{X}_i;h)+o_P(1)\stackrel{d}{\longrightarrow}\mathcal{N}(0,V^h_{tt'}),$$
where 
{\small
\begin{align}
&\psi_{tt'}(Y_i,T_i,\boldsymbol{M}_i,\boldsymbol{X}_i;h)=\frac{\pi_{\boldsymbol{M,X}}(T_i,\boldsymbol{M}_i,\boldsymbol{X}_i)\pi_{\boldsymbol{X}}(T_i+\delta,\boldsymbol{X}_i)}{\pi_{\boldsymbol{M,X}}(T_i+\delta,\boldsymbol{M}_i,\boldsymbol{X}_i)\cdot p_{t,h}}\mathcal{K}_h\left(T_i-t\right)\left\{Y_i-\mathbb{E}[Y|T_i,M_i,X_i] \right\}\notag,\\ 
&\quad\quad\quad\quad\quad -\frac{f_{T|\boldsymbol{X}}(T_i-\delta|\boldsymbol{X}_i)\pi_{\boldsymbol{X}}(T_i,\boldsymbol{X}_i)}{f_{T|\boldsymbol{X}}(T_i|\boldsymbol{X}_i)\cdot p_{t,h}}\mathcal{K}_h\left(T_i-t-\delta\right)\mathbb{E}\left[\frac{\pi_{\boldsymbol{M,X}}(T_i-\delta,\boldsymbol{M}_i,\boldsymbol{X}_i)Y_i}{\pi_{\boldsymbol{M,X}}(T_i,\boldsymbol{M}_i,\boldsymbol{X}_i)}\bigg|T_i,\boldsymbol{X}_i\right]\notag\\
&\quad\quad\quad\quad\quad +\frac{f_{T}(T_i-\delta)\pi_{\boldsymbol{X}}(T_i,\boldsymbol{X}_i)}{f_{T}(T_i)\cdot p_{t,h}}\mathcal{K}_h\left(T_i-t-\delta\right)\mathbb{E}\left[Y_i|T_i,\boldsymbol{M}_i,\boldsymbol{X}_i\right],\notag
\end{align}
and
\begin{align*}
V^h_{tt'}&=\lim_{h\to 0}h\cdot \mathbb{E}\left[\left\{\psi_{tt'}(Y_i,T_i,\boldsymbol{M}_i,\boldsymbol{X}_i;h)\right\}^2\right]\\
&=\frac{\kappa_{02}}{f_{T}(t)}\mathbb{E}\bigg[ \frac{\pi^2_{\boldsymbol{M,X}}(T_i,\boldsymbol{M}_i,\boldsymbol{X}_i)}{\pi^2_{\boldsymbol{M,X}}(T_i+\delta,\boldsymbol{M}_i,\boldsymbol{X}_i)}\pi^2_{\boldsymbol{X}}(T_i+\delta,\boldsymbol{X}_i)\left\{Y_i-\mathbb{E}[Y|T_i,\boldsymbol{M}_i,\boldsymbol{X}_i]\right\}^2\bigg| T_i=t\bigg]\notag\\
&\qquad+\frac{\kappa_{02}f_{T}(t+\delta)}{f^2_{T}(t)}\mathbb{E}\Bigg[\Bigg\{
\frac{\pi_{\boldsymbol{X}}(T_i,\boldsymbol{X}_i)f_{T|\boldsymbol{X}}(T_i-\delta|\boldsymbol{X}_i)}{f_{T|\boldsymbol{X}}(T_i|\boldsymbol{X}_i)}\mathbb{E}\left[\frac{\pi_{\boldsymbol{M,X}}(T_i-\delta,\boldsymbol{M}_i,\boldsymbol{X}_i)Y_i}{\pi_{\boldsymbol{M,X}}(T_i,\boldsymbol{M}_i,\boldsymbol{X}_i)} \bigg |T_i,\boldsymbol{X}_i\right]\notag\\
&\qquad\qquad\qquad-\frac{\pi_{\boldsymbol{X}}(T_i,\boldsymbol{X}_i)f_{T}(T_i-\delta)}{f_{T}(T_i)}\mathbb{E}\left[Y_i |T_i,\boldsymbol{M}_i,\boldsymbol{X}_i\right]\Bigg\}^2\bigg|T_i=t+\delta \Bigg],
\notag
\end{align*}}
and $p_{t,h}=\mathbb{E}\left[\mathcal{K}_h\left(T-t\right)\right]$, $\kappa_{ij}=\int u^i \mathcal{K}^j(u)du.$
\end{theorem}
Following Theorem \ref{thm:ker}, we establish the asymptotics of the estimated direct and indirect effect.
\begin{rk}\label{rk:EfficientvsHsu}
Note that the asymptotic variance of $\widehat{\mu}_{NKW}(t,t)$ and $\widehat{\mu}_{NKW}(t,t')$ introduced by \cite{hsu2018direct} are {\footnotesize
\begin{align*}
V^{NKW}_{tt}=\begin{cases}
\kappa_{02}\mathbb{E}\left[Var\left[Y|T=t,\boldsymbol{X}\right]/f_{T|\boldsymbol{X}}(t|\boldsymbol{X})\right] \ \ \ \text{if}\ h=h_1=h_2\ \text{and}\ \mathcal{K}_{1,h_1}(\cdot)=\mathcal{K}_{2,h_2}(\cdot);\\
\kappa_{02}\mathbb{E}\left\{\mathbb{E}\left[\{Y-\mu(t,t)\}^2|T=t,\boldsymbol{X} \right]/f_{T|\boldsymbol{X}}(t|\boldsymbol{X}) \right\}\ \ \text{if}\ h=h_2<h_1,
\end{cases}
\end{align*}}
and{\footnotesize
\begin{align*}
V^{NKW}_{tt'}=\begin{cases}
\kappa_{02}\Bigg(\mathbb{E}\left[Var(Y|T=t,\boldsymbol{M}, \boldsymbol{X})\frac{f_{T|\boldsymbol{M},\boldsymbol{X}}^2(t'|M,X)}{f_{T|\boldsymbol{M},\boldsymbol{X}}(t|\boldsymbol{M},\boldsymbol{X})f^2_{T|\boldsymbol{X}}(t'|\boldsymbol{X})} \right] \\
\qquad \quad+\mathbb{E}\left[Var\{g(t,\boldsymbol{M},\boldsymbol{X})|T=t',\boldsymbol{X}\}/f_{T|\boldsymbol{X}}(t'|\boldsymbol{X})\right]\Bigg)\ \ \text{if}\ h=h_1=h_2\ \text{and}\ \mathcal{K}_{1,h_1}(\cdot)=\mathcal{K}_{2,h_2}(\cdot);
\\
\kappa_{02}\mathbb{E}\left[\mathbb{E}\left\{(Y-\mu(t,t'))^2|T=t,\boldsymbol{M},\boldsymbol{X}\right\}\frac{f_{T|\boldsymbol{M},\boldsymbol{X}}^2(t'|\boldsymbol{M},\boldsymbol{X})}{f_{T|\boldsymbol{M},\boldsymbol{X}}(t|\boldsymbol{M},\boldsymbol{X})f^2_{T|\boldsymbol{X}}(t'|\boldsymbol{X})} \right]\ \ \text{if}\ h=h_2<h_1,
\end{cases}
\end{align*}}
where $\mathcal{K}_{1,h_1}$ and $\mathcal{K}_{2,h_2}$ are two prespecified kernel function used in the estimators of \cite{hsu2018direct}, $g(t,\boldsymbol{M},\boldsymbol{X})=\mathbb{E}[Y|T=t,\boldsymbol{M},\boldsymbol{X}]$.
We prove that $V^{h}_{tt}\le V_{tt}^{NKW}$ and $V^{h}_{tt'}\le V_{tt'}^{NKW}$ in section~S7 in the Supplementary Material.
\end{rk}
\end{document}